\documentclass{amsart}
\usepackage{amssymb}
\usepackage{mathrsfs}
\usepackage{stmaryrd}
\usepackage{bbm}
\usepackage{oldgerm}
\usepackage[english]{babel}
\usepackage[T1]{fontenc}
\usepackage[latin1]{inputenc}
\usepackage[all]{xy}
\usepackage{hyperref}
\usepackage{upref}
\usepackage{xcolor}
\usepackage{tikz-cd}

\hypersetup{
 colorlinks,
 linkcolor={red!50!black},
 citecolor={blue!50!black},
 urlcolor={blue!80!black}
}

\newtheorem{teo}[subsection]{Theorem}
\newtheorem{prop}[subsection]{Proposition}
\newtheorem{cor}[subsection]{Corollary}
\newtheorem{lem}[subsection]{Lemma}

\theoremstyle{definition}

\newtheorem{defi}[subsection]{Definition}
\newtheorem{rema}[subsection]{Remark}

\newtheorem{exemple}[subsection]{Example}

\numberwithin{equation}{subsection}

\newcommand{\A}{\mathcal{A}}

\newcommand{\Alg}{\mathrm{Alg}}

\newcommand{\Div}{\mathrm{Div}}
\newcommand{\E}{\mathcal{E}}

\newcommand{\et}{\mathrm{\acute{E}t}}
\newcommand{\eet}{\mathrm{\acute{e}t}}

\newcommand{\fp}{\mathrm{fp}}
\newcommand{\fppf}{\mathrm{fppf}}

\newcommand{\Fppf}{\mathrm{Fppf}}
\newcommand{\F}{\mathcal{F}}

\newcommand{\G}{\mathcal{G}}

\newcommand{\Hm}{\mathcal{H}}

\newcommand{\id}{\mathrm{id}}
\newcommand{\ka}{\kappa}
\newcommand{\Lc}{\mathcal{L}}
\newcommand{\M}{\mathcal{M}}

\newcommand{\Ow}{\mathcal{O}}

\newcommand{\Pic}{\mathrm{Pic}}

\newcommand{\Sch}{\mathrm{Sch}}

\newcommand{\Spec}{\mathrm{Spec}}
\newcommand{\Sy}{\mathrm{Sym}}
\newcommand{\Sym}{\mathrm{Sym}}

\newcommand{\Vc}{\mathcal{V}}

\newcommand{\Z}{\mathbb{Z}}

\title{On the ramified class field theory of relative curves}
\author{Quentin Guignard}
\address{Institut des
Hautes \'Etudes Scientifiques, 35 route de Chartres, 91440 Bures-sur-Yvette, France}
\address{
\'Ecole Normale Sup\'erieure, 45 rue d'Ulm, 75005 Paris, France}
\email{quentin.guignard@ens.fr}

\begin{document}

\begin{abstract}
We generalize Deligne's approach to tame geometric class field theory to the case of a relative curve, with arbitrary ramification. 
\end{abstract}

\maketitle
\tableofcontents

\section{Introduction}

Let $X \rightarrow S$ be a relative curve, i.e. a smooth morphism of schemes of relative dimension $1$, with connected geometric fibers, which is Zariski-locally projective over $S$. Let $Y \hookrightarrow X$ be a relative effective Cartier divisor over $S$ (cf. \ref{chap21.2}), and let $U$ be the complement of $Y$ in $X$.

The pairs $(\Lc, \alpha)$, where $\Lc$ is an invertible $\Ow_X$-module and $\alpha$ is a rigidification of $\Lc$ along $Y$, are parametrized by an $S$-group scheme $\Pic_S(X,Y)$, the relative rigidified Picard scheme (cf. \ref{chap2rigidrepr}). The Abel-Jacobi morphism
$$
\Phi : U \rightarrow \Pic_S(X,Y)
$$
is the morphism which sends a section $x$ of $U$ to the pair $(\Ow(x),1)$, cf. \ref{chap2abeljacobi}. We prove the following relative version of the main theorem of geometric global class field theory:

\begin{teo}\label{chap2teoteo}(Th. \ref{chap2GCFT}) Let $\Lambda$ be a finite ring of cardinality invertible on $S$, and let $\F$ be an \'etale sheaf of $\Lambda$-modules, locally free of rank $1$ on $U$, with ramification bounded by $Y$ (cf. \ref{chap2ramif6}). Then, there exists a unique (up to isomorphism) multiplicative \'etale sheaf of $\Lambda$-modules $\G$ on $\Pic_S(X,Y)$, locally free of rank $1$, such that the pullback of $\G$ by $\Phi$ is isomorphic to $\F$.
\end{teo}

The notion of multiplicative locally free $\Lambda$-module of rank $1$ is defined in \ref{chap2chartors}, and it corresponds to isogenies $G \rightarrow \Pic_S(X,Y)$ with constant kernel $\Lambda^{\times}$. We restrict ourself in this article to $\Lambda^{\times}$-torsors, with $\Lambda$ as in Theorem \ref{chap2teoteo}, in order to simplify the exposition, since we are able to apply directly our main descent tool in this context, namely Lemma \ref{chap2pull}. However, the latter lemma, and hence Theorem \ref{chap2teoteo} can be extended to $G$-torsors, where $G$ is an arbitrary locally constant finite abelian group on $S_{\text{\'et}}$.

The case where $S$ is the spectrum of a perfect field is originally due Serre and Lang, cf. (\cite{Lang}, 6) and \cite{JPS3}. Their proof relies on the Albanese property of Rosenlicht's generalized Jacobians \cite{Ros}. A similar proof was sketched by Deligne in his $1974$ letter to Serre, published as an appendix in \cite{bloch}.  However, a more geometric proof was given by Deligne in the tamely ramified case; an account of Deligne's proof in the unramified case over a finite field can be found in (\cite{Laumon}, Sect. $2$). We generalize the latter approach by Deligne to allow arbitrary ramification and an arbitrary base $S$. This generalization is inspired by notes by Alain Genestier (unpublished) on arithmetic global class field theory.

Deligne's approach has the advantage over Serre and Lang's to yield an explicit geometric construction of the isogeny over $\Pic_S(X,Y)$ corresponding to a local system of rank $1$ over $U$. This feature of Deligne's approach carries over to ours, and is in fact crucial in order to handle the case of an arbitrary base $S$. 

The author was informed during the preparation of this manuscript that Daichi Takeuchi has independently obtained a different proof of \ref{chap2teoteo} in the case where $S$ is the spectrum of a perfect field, also by generalizing Deligne's approach to handle arbitrary ramification, cf. \cite{T18}.

\subsection*{Acknowledgements} This work is part of the author's PhD dissertation and it was prepared at the Institut des Hautes \'Etudes Scientifiques and the \'Ecole Normale Sup\'erieure while the author benefited from their hospitality and support. The author is indebted to Ahmed Abbes for his numerous suggestions and remarks. Further thanks go to Salim Tayou and Christophe Soul\'e for their comments on this text. The author would also like to thank Fu Lei, Fabrice Orgogozo and the anonymous referee for remarks and suggestions which helped to improve this text.

\subsection*{Notation and conventions.} We fix a universe $\mathcal{U}$ (\cite{SGA4}, I.0). Thoughout this paper, all sets are assumed to belong to $\mathcal{U}$ and we will use the term ``topos'' as a shorthand for ``$\mathcal{U}$-topos'' (\cite{SGA4}, IV.1.1). The category of sets belonging to $\mathcal{U}$ is simply denoted by $\mathrm{Sets}$.

For any integers $e,d$ we denote by $\llbracket e,d\rrbracket$ the set of integers $n$ such that $e \leq n \leq d$ and by $\mathfrak{S}_d$ the group of bijections of $\llbracket1,d\rrbracket$ onto itself. 

In this paper, all rings are unital and commutative. For any ring $A$, we denote by $\Alg_A$ the category of $A$-algebras. For any scheme $S$, we denote by $\Sch_{/S}$ the category of $S$-schemes. We denote by $S_{\text{\'et}}$ (resp. $S_{\et}$) the small \'etale topos (resp. big \'etale topos) of a scheme $S$, i.e. the topos of sheaves of sets for the \'etale topology (\cite{SGA4}, VII.1.2) on the category of \'etale $S$-schemes (resp. on $\Sch_{/S}$), and by $S_{\Fppf}$ the big $\fppf$ topos of $S$, i.e. the topos of sheaves of sets for the $\fppf$ topology on $\Sch_{/S}$ (\cite{SGA4}, VII.4.2). If $f : X \rightarrow S$ is a morphism of schemes, then we denote by $(f^{-1},f_*)$ the induced morphism of topos from $X_{\et}$ to $S_{\et}$. The symbol $f^*$ will exclusively denote the pullback functor from $\Ow_S$-modules to $\Ow_X$-modules.

\section{Preliminaries}

\subsection{\label{chap23.0}}Let $E$ be a topos and let $G$ be an abelian group in $E$. We denote by $G E$ the category of objects of $E$ endowed with a left action of $G$. If $X$ is an object of $E$, we denote by $E_{/X}$ the topos of $X$-objects in $E$. If $X$ is considered as an object of $G E$ by endowing it with the trivial left $G$-action, then we have $(G E)_{/X} = G(E_{/X})$ and this category will be simply denoted by $G E_{/X}$.

\begin{defi}\label{chap2defitorseur} A \textbf{$G$-torsor} over an object $X$ of $E$ is an object $P$ of $G E_{/X}$ such that $P \rightarrow X$ is an epimorphism and the morphism
\begin{align*}
G\times_X P &\rightarrow P \times_X P \\
(g,p) &\mapsto (g \cdot p,p)
\end{align*}
is an isomorphism. We denote by $\mathrm{Tors}(X,G)$ the full subcategory of $G E_{/X}$ whose objects are the $G$-torsors over $X$. If $f : Y \rightarrow X$ is a morphism in $E$, we denote by $f^{-1} : \mathrm{Tors}(X,G) \rightarrow \mathrm{Tors}(Y,G) $ the functor which associates $f^{-1} P = P \times_{X,f} Y$ to a $G$-torsor $P$ over $X$. 
\end{defi}

The category $\mathrm{Tors}(X,G)$ is monoidal, with product
$$
P_1 \otimes P_2 = G_2 \setminus P_1 \times_X P_2,
$$
where $G_2$ is the kernel of the multiplication morphism $G \times G \rightarrow G$, and where $G_2 \hookrightarrow G \times G$ acts diagonally on $P_1 \times_X P_2$. The neutral element for this product is the trivial $G$-torsor over $X$, namely $ G \times X$, and each $G$-torsor $P$ over $X$ is invertible with respect to $\otimes$, with inverse given by
$$
P^{-1} = \underline{\mathrm{Hom}}_{G E_{/X}}(P, G \times X),
$$
where $\underline{\mathrm{Hom}}_{G E_{/X}}$ denotes the internal $\mathrm{Hom}$ functor in $G E_{/X}$.

\begin{exemple}\label{chap2ex1} If $G = \Lambda^{\times}$ for some ring $\Lambda$ in $E$, then the monoidal category $\mathrm{Tors}(X,G)$ is equivalent to the groupoid of locally free $\Lambda$-modules of rank $1$ in $E_{/X}$. The equivalence is given by the functor which sends an object $P$ of $\mathrm{Tors}(X,G)$ to the $\Lambda$-module $G \setminus (\Lambda \times P)$, where the action of $G = \Lambda^{\times}$ on $\Lambda \times P$ is given by the formula $g \cdot (\lambda,p) = (g \lambda, g \cdot p)$. The functor which sends a locally free $\Lambda$-module $M$ of rank $1$ of $E_{/X}$ to the $G$-torsor of isomorphisms of $\Lambda$-modules from $M$ to $\Lambda$ defines a quasi-inverse to the latter functor.
\end{exemple}

\subsection{\label{chap23.1}} Let $E$ be a topos, and let us denote by $1$ its terminal object. Let us consider an exact sequence
$$
1 \rightarrow G \xrightarrow[]{i} P \xrightarrow[]{r} Q \rightarrow 1
$$
of abelian groups in $E$. The morphism 
\begin{align*}
G \times_Q P &\rightarrow P \times_Q P \\
(g, p) &\mapsto (i(g) + p,p)
\end{align*}
is an isomorphism, so that $P$ is a $G$-torsor over $Q$. Moreover, the multiplication morphism
$$
P \times P \rightarrow P
$$
factors though $G_2 \setminus P \times P $, where $G_2 \hookrightarrow G \times G$ is the kernel of the multiplication morphism of $G$, acting diagonally on $P \times P$. We thus obtain a morphism
$$
p_1^{-1} P \otimes p_2^{-1} P \rightarrow m^{-1} P
$$
of $G$-torsors over $Q \times Q$, where $p_1$ and $p_2$ are the canonical projections and $m$ is the multiplication morphism of $Q$. 

The following definition is inspired by (\cite{MB}, I.2.3):

\begin{defi}\label{chap2chartors} Let $G$ be an abelian group of $E$ and let $Q$ be a commutative semigroup of $E$ (with or without identity). Let $m : Q \times Q \rightarrow Q$ be the multiplication morphism of $Q$. A \textbf{multiplicative $G$-torsor} over $Q$ is a $G$-torsor $P \rightarrow Q$, together with an isomorphism $\theta : p_1^{-1} P \otimes p_2^{-1} P \rightarrow m^{-1} P$ of $G$-torsors over $Q \times Q$ where $p_1$ and $p_2$ are the canonical projections, which satisfy the following two properties.
\begin{itemize}
\item[$\triangleright$] \textbf{Symmetry:} if $\sigma$ is the involution of $Q \times Q$ which switches the two factors, then the isomorphism
$$
 p_2^{-1} P \otimes p_1^{-1} P \rightarrow \sigma^{-1}(p_1^{-1} P \otimes p_2^{-1} P) \xrightarrow[]{\sigma^{-1} \theta} \sigma^{-1} m^{-1} P \rightarrow m^{-1} P
$$
is the composition of $\theta$ with the canonical isomorphism $ p_2^{-1} P \otimes p_1^{-1} P \rightarrow p_1^{-1} P \otimes p_2^{-1} P$.
\item[$\triangleright$] \textbf{Associativity:} if $q_i : Q \times Q \times Q \rightarrow Q$ (resp. $q_{ij} : Q \times Q \times Q \rightarrow Q \times Q$) is the projection on the $i$-th factor for $i \in \llbracket1,3\rrbracket$ (resp. on the $i$-th and $j$-th factors for $(i,j) \in \llbracket1,3\rrbracket^2$ such that $i<j$) and if $m_3 : Q \times Q \times Q \rightarrow Q$ is the multiplication morphism, then the diagram of $G$-torsors over $Q \times Q \times Q$
\begin{center}
 \begin{tikzpicture}[scale=1]

\node (A) at (0,2) {$q_1^{-1} P \otimes q_2^{-1} P \otimes q_3^{-1} P$};
\node (B) at (3,4) {$q_1^{-1} P \otimes (m q_{23})^{-1} P$};
\node (C) at (3,0) {$(m q_{12})^{-1} P \otimes q_3^{-1} P$};
\node (D) at (6,2) {$m_3^{-1} P$};
\path[->,font=\scriptsize]
(A) edge node[above left]{$id \otimes q_{23}^{-1} \theta $} (B)
(A) edge node[below left]{$ q_{12}^{-1} \theta \otimes  id$} (C)
(B) edge node[above right]{$ (q_1 \times mq_{23})^{-1} \theta $} (D)
(C) edge node[below right]{$ (mq_{12} \times q_3)^{-1} \theta$} (D);

\end{tikzpicture} 
\end{center}
is commutative.
\end{itemize}
The category of multiplicative $G$-torsors is fibered in groupoids over the category of commutative semigroups of $E$. We denote by $\mathrm{Tors}^{\otimes}(Q,G)$ the groupoid of multiplicative $G$-torsors over a commutative semigroup $Q$ of $E$.
\end{defi}

\begin{rema}\label{chap2charsheaf} If $G = \Lambda^{\times}$ for some ring $\Lambda$ in $E$, we use the term  ``\textbf{multiplicative locally free $\Lambda$-module of rank $1$}'' as a synonym for ``multiplicative $G$-torsor'', when we want to emphasize the locally free $\Lambda$-module of rank $1$ corresponding to a given $G$-torsor, rather than the $G$-torsor itself (cf. \ref{chap2ex1}). 
\end{rema}

\begin{prop}\label{chap2mondesc2}Let $G$ be an abelian group in $E$, let $Q$ be a commutative semigroup in $E$ and let $I$ be an ideal of $Q$. If the projection morphisms $Q \times I \rightarrow Q$ and $I \times I \rightarrow I$ onto the first factors are morphisms of descent for the fibered category of multiplicative $G$-torsors (cf. \ref{chap2chartors}), then the restriction functor
$$
\mathrm{Tors}^{\otimes}(Q,G) \rightarrow \mathrm{Tors}^{\otimes}(I,G)
$$
is fully faithful.
\end{prop}

Let $i : I \rightarrow Q$ be the canonical injection morphism. Let $p_1$ and $p_2$ be the projection morphisms of $Q \times I$ onto its first and second factors respectively, and let $m : Q \times I \rightarrow I$ be the multiplication morphism. Let $(P,\theta)$ and $(P',\theta')$ be multiplicative $G$-torsors over $Q$. We have an isomorphism
$$
\beta_P : p_1^{-1} P \xrightarrow[]{(\mathrm{id} \times i)^{-1} \theta} m^{-1} i^{-1}P \otimes p_2^{-1} i^{-1}P^{-1},
$$
and similarly for $P'$. If $\alpha : i^{-1}P \rightarrow i^{-1}P'$ is a morphism of multiplicative $G$-torsors over $I$, then $\beta_{P'}^{-1}(m^{-1}\alpha \otimes p_2^{-1} \alpha) \beta_P$ is an isomorphism from $p_1^{-1} P$ to $p_1^{-1} P'$, which is compatible with the canonical descent datum for $p_1$ associated to $p_1^{-1} P$ and $p_1^{-1} P'$: indeed, if $q_1 : Q \times I \times I \rightarrow Q$ and $q_2,q_3 : Q \times I \times I \rightarrow I$ (resp. $q_{ij}$) are the projections on the first, second and third factors of $ Q \times I \times I$ respectively (resp. on the product of its $i$-th and $j$-th factors for $(i,j) \in \llbracket1,3\rrbracket^2$ such that $i<j$) and if $m_3 : Q \times I \times I \rightarrow I$ is the multiplication morphism, then the diagram of $G$-torsors over $Q \times I \times I$
\begin{center}
 \begin{tikzpicture}[scale=1]

\node (A) at (0,2) {$q_{12}^{-1} p_1^{-1} P$};
\node (B) at (4,4) {$(im_3)^{-1} P \otimes (im q_{23})^{-1} P^{-1}$};
\node (C) at (4,0) {$(imq_{12})^{-1} P \otimes (ip_2q_{12})^{-1}  P^{-1}$};
\node (D) at (8,2) {$(imq_{12})^{-1} P \otimes (iq_3)^{-1}  P \otimes (ip_2q_{12})^{-1} P^{-1}\otimes (iq_3)^{-1} P^{-1}$};
\path[->,font=\scriptsize]
(A) edge node[above left]{$(q_1 \times mq_{23})^{-1} \theta $} (B)
(A) edge node[below left]{$ q_{12}^{-1} \beta_P$} (C)
(B) edge node[above right]{$ (mq_{12} \times q_3)^{-1} \theta^{-1} \otimes q_{23}^{-1} \theta^{-1}  $} (D)
(C) edge  (D);

\end{tikzpicture} 
\end{center}
is commutative, and similarly for $P'$, so that the pullback of $\beta_{P'}^{-1}(m^{-1}\alpha \otimes p_2^{-1} \alpha) \beta_P$ by $q_{12}$ is given by the composition
$$
q_{1}^{-1}  P \rightarrow (im_3)^{-1}  P \otimes (im q_{23})^{-1}  P^{-1}  \xrightarrow[]{m_3^{-1} \alpha \otimes (m q_{23})^{-1} \alpha}   (im_3)^{-1}P\otimes (im q_{23})^{-1} P'^{\prime-1}  \rightarrow  q_{1}^{-1}  P'  ,
$$
and therefore coincides with its pullback by $q_{13}$.

 Since $p_1$ is a morphism of descent for the fibered category of multiplicative $G$-torsors, there is a unique morphism $\gamma : P \rightarrow P'$ of multiplicative $G$-torsors over $Q$ such that $p_1^{-1}\gamma = \beta_{P'}^{-1} (m^{-1}\alpha \otimes p_2^{-1} \alpha) \beta_P$. The restriction of $p_1^{-1} \gamma$ to $I \times I$ is the pullback of $\alpha$ by the first projection, which is a morphism of descent for the fibered category of multiplicative $G$-torsors, so that the restriction of $\gamma$ to $I$ is $\alpha$.

\begin{prop}\label{chap2mondesc}Let $G$ be an abelian group in $E$, and let $\rho : M \rightarrow Q$ be a morphism of commutative semigroups in $E$. If $\rho$ (resp. $\rho \times \rho$ and $\rho  \times \rho \times \rho$) is a morphism of effective descent (resp. of descent) for the fibered category of $G$-torsors, then $\rho$ is a morphism of effective descent for the fibered category of multiplicative $G$-torsors.
\end{prop}

A descent datum of multiplicative $G$-torsors for $\rho$ yields a descent datum of $G$-torsors for $\rho$, hence a $G$-torsor over $Q$ by hypothesis. Since $\rho \times \rho$ and $\rho  \times \rho \times \rho$ are morphisms of descent for the fibered category of $G$-torsors, the structure of multiplicative $G$-torsor descends as well. Details are omitted.

\begin{prop}\label{chap2ext} Let $G$ and $Q$ be abelian groups in $E$. The groupoid $\mathrm{Tors}^{\otimes}(Q,G)$ of multiplicative $G$-torsors over $Q$ is equivalent as a monoidal category to the groupoid of extensions of $Q$ by $G$ in $E$, with the Baer sum as a monoidal structure.
\end{prop}

We have already seen how to associate a multiplicative $G$-torsor to an extension of $Q$ by $G$. This construction is functorial, and the corresponding functor is an equivalence by (\cite{MB}, I.2.3.10).

\begin{cor}\label{chap2ext2} Let $G$ and $Q$ be abelian groups in $E$. The group of isomorphism classes of multiplicative $G$-torsors over $Q$ is isomorphic to the group $\mathrm{Ext}^1(Q,G)$ of isomorphism classes of extensions of $Q$ by $G$ in $E$.
\end{cor}

\subsection{\label{chap23.14}}Let $S$ be a scheme, let $X$ be an $S$-scheme, and let $G$ be a finite abelian group. Let $P$ be a $G$-torsor over $X$ in $S_{\et}$. Since $P \rightarrow X$ is an epimorphism in $S_{\et}$, there is an \'etale cover $(X_i \rightarrow X)_{i \in I}$ such that for each $i \in I$, the morphism $X_i \rightarrow X$ factors through $P \rightarrow X$. In particular, for each $i \in I$ the $G$-torsor $P \times_X X_i \rightarrow X_i$ is isomorphic to the trivial $G$-torsor $G \times X_i \rightarrow X_i$, so that $P \times_X X_i $ is representable by a finite \'etale $X_i$-scheme. By  \'etale descent of affine morphisms, we obtain:

\begin{prop}\label{chap2torsrepre} Let $G$ be a finite abelian group, let $S$ be a scheme, and let $P$ be a $G$-torsor over an $S$-scheme $X$ in $S_{\et}$. Then the \'etale sheaf $P : \mathrm{Sch}_{/S} \rightarrow \mathrm{Sets}$ is representable by a finite \'etale $X$-scheme.
\end{prop}

The topos $(S_{\et})_{/X}$ coincides with $X_{\et}$. The category of $G$-torsors over $X$ in $S_{\et}$ is therefore equivalent to the category of $G$-torsors over the terminal object in $X_{\et}$, and Proposition \ref{chap2torsrepre} yields:

\begin{cor}\label{chap22.13} Let $G$ be a finite abelian group, let $S$ be a scheme, and let $X$ be an $S$-scheme. Then the category of $G$-torsors over $X$ in $S_{\et}$ is equivalent to the category of $G$-torsors over the terminal object in $X_{\text{\'et}}$.
\end{cor}

\subsection{\label{chap23.15}}Let $S$ be a scheme, and let $G$ be a finite abelian group. Let $Q$ be a commutative $S$-group scheme, and let $M$ be a sub-$S$-semigroup scheme of $Q$. 

\begin{prop}\label{chap2extension} Assume that the morphism
\begin{align*}
\rho : M \times_S M &\rightarrow Q \\
(x,y) &\mapsto x y^{-1}
\end{align*}
is faithfully flat and quasi-compact, and that $M$ is flat over $S$. Then the restriction functor
$$
\mathrm{Tors}^{\otimes}(Q,G) \rightarrow \mathrm{Tors}^{\otimes}(M,G),
$$
is an equivalence of categories.
\end{prop}

Let $(P,\theta)$ be a multiplicative $G$-torsor over $M$. For $i \in \llbracket1,4\rrbracket$, let $r_i$ be the projection of $R = (M \times_S M) \times_{\rho,Q,\rho} (M \times_S M)$ onto its $i$-th factor. Similarly, for $i,j \in \llbracket1,4\rrbracket$ such that $i<j$, we denote by $r_{ij} : R \rightarrow M \times_S M$ the projection on the $i$-th and $j$-th factors. We then have a sequence of isomorphisms
\begin{align*}
\left( r_1^{-1} P \otimes r_2^{-1} P^{-1} \right) \otimes \left( r_3^{-1} P \otimes r_4^{-1} P^{-1} \right)^{-1} &\rightarrow r_{14}^{-1}(p_1^{-1} P \otimes p_2^{-1} P) \otimes r_{23}^{-1}(p_1^{-1} P \otimes p_2^{-1} P)^{-1} \\
&\xrightarrow[]{r_{14}^{-1} \theta \otimes (r_{23}^{-1} \theta )^{-1}}  (m r_{14})^{-1} P   \otimes ((m r_{23})^{-1} P)^{-1},
\end{align*}
of $G$-torsors over $R$, where $m : M \times_S M \rightarrow M$ is the multiplication of $M$. Since $m r_{14} = m r_{23}$, the latter $G$-torsor is canonically trivial. We thus obtain an isomorphism
$$
\psi :  r_1^{-1} P \otimes r_2^{-1} P^{-1} \rightarrow r_3^{-1} P \otimes r_4^{-1} P^{-1},
$$
of $G$-torsors over $R$. The associativity of $\theta$ (cf. \ref{chap2chartors}) implies that $\psi$ is a cocycle, i.e. $(p_1^{-1} P \otimes p_2^{-1} P^{-1}, \psi)$ is a descent datum for $\rho$. By Proposition \ref{chap2torsrepre} and since faithfully flat and quasi-compact morphisms of schemes are of effective descent for the fibered category of affine morphisms, the conditions of Proposition \ref{chap2mondesc} are satisfied, and thus there exists a multiplicative $G$-torsor $P'$ over $Q$ and an isomorphism $\alpha : \rho^{-1} P' \rightarrow p_1^{-1} P \otimes p_2^{-1} P^{-1}$ such that $\psi$ is given by the composition
$$
 r_1^{-1} P \otimes r_2^{-1} P^{-1} \xrightarrow[]{r_{12}^{-1} \alpha^{-1}} (\rho r_{12})^{-1} P' = (\rho r_{34})^{-1} P'   \xrightarrow[]{r_{34}^{-1} \alpha}  r_3^{-1} P \otimes r_4^{-1} P^{-1}.
$$
The association $P \mapsto P'$ then defines a functor from $ \mathrm{Tors}^{\otimes}(M,G)$ to $ \mathrm{Tors}^{\otimes}(Q,G)$. For any multiplicative $G$-torsor $U$ over $Q$, we have an isomorphism $U \rightarrow (U\times_Q M)'$ by multiplicativity, which is functorial in $U$.

We now construct, for any multiplicative $G$-torsor $(P,\theta)$ over $M$, an isomorphism $P \rightarrow P' \times_Q M$ of multiplicative $G$-torsors which is functorial in $P$. Let $\nu : M \times_S M \rightarrow M \times_S M$ be the morphism which sends a section $(x,y)$ to $(xy,y)$. We have an isomorphism
$$
(\rho \nu)^{-1} P' \xrightarrow[]{\nu^{-1} \alpha} \nu^{-1}(p_1^{-1} P \otimes p_2^{-1} P^{-1}) \rightarrow m^{-1} P \otimes p_2^{-1} P^{-1} \xrightarrow[]{\theta^{-1}} p_1^{-1} P.
$$
The diagram 
\begin{center}
 \begin{tikzpicture}[scale=1]

\node (A) at (0,1) {$M \times_S M$};
\node (B) at (3,2) {$M \times_S M$};
\node (C) at (3,0) {$M$};
\node (D) at (6,1) {$Q$};
\path[->,font=\scriptsize]
(A) edge node[above left]{$\nu$} (B)
(A) edge node[below left]{$p_1$} (C)
(B) edge node[above right]{$ \rho $} (D)
(C) edge node[below right]{} (D);

\end{tikzpicture} 
\end{center}
is commutative, hence $(\rho \nu)^{-1} P'$ is isomorphic to $p_1^{-1}(P' \times_Q M)$. We thus obtain an isomorphism
$$
\beta : p_1^{-1} P \rightarrow p_1^{-1}(P' \times_Q M),
$$
of multiplicative $G$-torsors. The morphism $\beta$ is compatible with the canonical descent data for $p_1$ associated to $p_1^{-1} P$ and $p_1^{-1}(P' \times_Q M)$. Since $p_1$ is a covering for the $\mathrm{fpqc}$ topology, Proposition \ref{chap2mondesc} applies, hence there is a unique isomorphism $\gamma : P \rightarrow P' \times_Q M$ of multiplicative $G$-torsors such that $\beta = p_1^{-1} \gamma$. The construction of this isomorphism of multiplicative $G$-torsors is functorial in $P$, hence the result.

\subsection{\label{chap23.10}}Let $A$ be a ring. If $M$ is an $A$-module, we denote by $\underline{M}$ the functor $B \mapsto M \otimes_A B$ from $\mathrm{Alg}_A$ to $\mathrm{Sets}$.

\begin{defi}(\cite{SGA4}, XVII 5.5.2.2)\label{chap2homoge} Let $M$ and $N$ be $A$-modules. A \textbf{polynomial map} from $M$ to $N$ is a morphism of functors $\underline{M} \rightarrow \underline{N}$. A polynomial map $f : \underline{M} \rightarrow \underline{N}$ is \textbf{homogeneous of degree $d$} if for any $A$-algebra $B$, any element $\lambda$ of $B$ and any element $m$ of $\underline{M}(B)$, we have $f(\lambda m) = \lambda^d f(m)$.
\end{defi}

For each integer $d$ and any $A$-module $M$, let $\mathrm{TS}_A^d(M) = (M^{\otimes_A d})^{\mathfrak{S}_d}$ be the $A$-module of symmetric tensors of degree $d$ with coefficients in $M$. If $M$ is a free $A$-module with basis $(e_i)_{i \in I}$, then we have a decomposition
\begin{align}\label{chap2dec}
TS_A^d(M) = \left( \bigoplus_{\beta : \llbracket 1,d\rrbracket \rightarrow I } A e_{\beta(1)} \otimes \dots \otimes e_{\beta(d)} \right)^{\mathfrak{S}_d} =  \bigoplus_{\substack{\alpha : I \rightarrow \mathbb{N} \\ \sum_{i \in I} \alpha(i) = d }} A e_{\alpha},
\end{align}
where we have set
$$
e_{\alpha} = \sum_{\substack{\beta : \llbracket 1,d\rrbracket \rightarrow I \\ \forall i, |\beta^{-1}(\{i\})| = \alpha(i) }} e_{\beta(1)} \otimes \dots \otimes e_{\beta(d)}.
$$
In particular $TS_A^d(M)$ is a free $A$-module, and its formation commutes with base change by any ring morphism $A \rightarrow B$.

\begin{prop}\label{chap2flat} Let $M$ be a flat $A$-module and let $d \geq 0$ be an integer. Then $\mathrm{TS}_A^d(M)$ is a flat module, and for any $A$-algebra $B$ the canonical homomorphism
$$
\mathrm{TS}_A^d(M) \otimes_A B \rightarrow \mathrm{TS}_B^d(M \otimes_A B)
$$
is bijective.
\end{prop}

Any flat $A$-module is a filtered colimit of finite free modules. We have already seen that the conclusion of Proposition \ref{chap2flat} holds whenever $M$ is free, hence the conclusion in general since the functor $\mathrm{TS}_A^d$ commutes with filtered colimits.

\begin{prop}\label{chap2flat2} Let $M$ be a flat $A$-module and let $d \geq 0$ be an integer. Let $\gamma_{d} : \underline{M} \rightarrow \underline{\mathrm{TS}_A^d(M)}$ be the functor which sends, for any $A$-algebra $B$, an element $m$ of $\underline{M}(B)$ to the element $m^{\otimes d}$ of $\mathrm{TS}_B^d(M \otimes_A B) = \mathrm{TS}_A^d(M) \otimes_A B$ (cf. \ref{chap2flat}). Then, for any homogeneous polynomial map $f : \underline{M} \rightarrow \underline{N}$ of degree $d$, there is a unique $A$-linear homomorphism $\tilde{f} : \mathrm{TS}_A^d(M) \rightarrow N$ such that $f = \tilde{f} \gamma_d$.
\end{prop}

As in Proposition \ref{chap2flat}, we can assume that $M$ is free of finite rank over $A$. Let $(e_i)_{i \in I}$ be a basis of $M$. Let us write
$$
f \left( \sum_{i \in I} X_i e_i \right) = \sum_{\substack{\alpha : I \rightarrow \mathbb{N} \\  }} X^{\alpha} f_{\alpha} 
$$
in $\underline{N}(A[(X_i)_{i \in I}])$ for some elements $(f_{\alpha})_{\alpha}$ of $N$, where $X^{\alpha} = \prod_{i \in I} X_i^{\alpha_i}$. Accordingly, we have for any $A$-algebra $B$ and any element $m = \sum_{i \in I} b_i e_i$ of $\underline{M}(B)$, the formula
$$
f \left( m \right) = \sum_{\substack{\alpha : I \rightarrow \mathbb{N} \\  }}  b^{\alpha} f_{\alpha},
$$
where $b^{\alpha} = \prod_{i \in I} b_i^{\alpha_i}$, by using the naturality of $f$ with the unique morphism of $A$-algebras $A[(X_i)_{i \in I}] \rightarrow B$ which sends $X_i$ to $b_i$ for each $i$. By applying this to the element $m = \sum_{i \in I} T X_i e_i$ of $\underline{M}(A[T,(X_i)_{i \in I}])$, we obtain
$$
f \left(\sum_{i \in I} T X_i e_i \right) = \sum_{\substack{\alpha : I \rightarrow \mathbb{N} \\ }} T^{|\alpha|} X^{\alpha} f_{\alpha},
$$
where we have set $|\alpha| = \sum_{i \in I} \alpha(i)$. Since $f$ is homogeneous of degree $d$, the left side of this equation is also equal to
$$
T^d f \left( \sum_{i \in I} X_i e_i \right) = \sum_{\substack{\alpha : I \rightarrow \mathbb{N} \\  }} T^d X^{\alpha} f_{\alpha}.
$$
We conclude that $T^d f_{\alpha} = T^{|\alpha|} f_{\alpha}$ in $N \otimes_A A[T]$ for any $\alpha : I \rightarrow \mathbb{N}$, and thus that $f_{\alpha}= 0$ whenever $|\alpha|$ differs from $d$. We therefore have
$$
f \left( m \right) = \sum_{\substack{\alpha : I \rightarrow \mathbb{N} \\ |\alpha|= d }}  b^{\alpha} f_{\alpha},
$$
for any $A$-algebra $B$ and any element $m = \sum_{i \in I} b_i e_i$ of $\underline{M}(B)$. Using the decomposition (\ref{chap2dec}), we also have
$$
\gamma_d (m) = \sum_{\beta : \llbracket 1,d\rrbracket \rightarrow I } \otimes_{j=1}^d b_{\beta(j)} e_{\beta(j)} = \sum_{\substack{\alpha : I \rightarrow \mathbb{N} \\ |\alpha| = d }}  b^{\alpha} e_{\alpha}.
$$ 
The conclusion of Proposition \ref{chap2flat2} is achieved by taking $\tilde{f}$ to be the unique morphism of $A$-modules from $\mathrm{TS}_A^d(M)$ to $N$ which sends $e_{\alpha}$ to $f_{\alpha}$.

\subsection{\label{chap23.11}}Let $A \rightarrow C$ be a ring morphism such that $C$ is a finitely generated projective $A$-module of rank $d$. For any $A$-algebra $B$ and any element $c$ of $\underline{C}(B)$, we set
$$
N_{C/A}(c) =  \mathrm{det}_{\underline{A}(B)}(m_c),
$$
where $m_c$ is the $\underline{A}(B)$-linear endomorphism of $ \underline{C}(B)$ induced by the multiplication by $c$. This defines a homogeneous polynomial map $N_{C/A} : \underline{C} \rightarrow \underline{A}$ of degree $d$ (cf. \ref{chap2homoge}). By \ref{chap2flat2}, there is a unique morphism of $A$-modules $\varphi : \mathrm{TS}_A^d(C) \rightarrow A$ such that $N_{C/A} = \varphi \gamma_d$.

\begin{prop}[\cite{SGA4}, XVII 6.3.1.6]\label{chap2morphalg} The morphism of $A$-modules $\varphi : \mathrm{TS}_A^d(C) \rightarrow A$ is a morphism of $A$-algebras.
\end{prop}

Let $x$ be an element of $C$, and let us consider the morphism of $A$-modules $f : y  \rightarrow \varphi(\gamma_d(x)y)$ from $\mathrm{TS}_A^d(C)$ to $A$. For any $A$-algebra $B$ and any element $c$ of $\underline{C}(B)$, we have
$$
f(\gamma_d(c)) =  \varphi(\gamma_d(x)\gamma_d(c)) = \varphi(\gamma_d(xc)) = N_{C/A}(xc) = N_{C/A}(x) N_{C/A}(c)
$$
by the multiplicativity of determinants, so that $f(\gamma_d(c)) = N_{C/A}(x) \varphi(\gamma_d(c))$. By the uniqueness statement in \ref{chap2flat2}, we obtain $f = N_{C/A}(x) \varphi$, i.e. for all $y$ in $\mathrm{TS}_A^d(C)$ we have
\begin{align}\label{chap2prelim}
\varphi(\gamma_d(x)y) = N_{C/A}(x) \varphi(y).
\end{align}
For any $A$-algebra $B$, one can apply this argument to the morphism $B \rightarrow \underline{C}(B)$ instead of $A \rightarrow C$. Thus (\ref{chap2prelim}) also holds for any element $x$ of $\underline{C}(B)$ and any element $y$ of $\underline{\mathrm{TS}_A^d(C)}(B) = \mathrm{TS}_{\underline{A}(B)}^d(\underline{C}(B))$ (cf. \ref{chap2flat}).
Now, let $y$ be an element of $\mathrm{TS}_A^d(C)$ and let us consider the morphism of $A$-modules $g : z \rightarrow \varphi(zy)$ from $ \mathrm{TS}_A^d(C) $ to $A$. We have proved that $g \gamma_d = \varphi(y) N_{C/A}$, hence $g = \varphi(y) \varphi$ by \ref{chap2flat2}. Thus $\varphi$ is a morphism of rings. Since $\varphi$ is also $A$-linear, it is a morphism of $A$-algebras.

\subsection{\label{chap23.4}} Let $S$ be a scheme.

\begin{defi}[\cite{SGA1}, V.1.7]\label{chap2defquot}
\end{defi}

\begin{itemize}

\item[$\triangleright$]Let $T$ be an object of a category $C$ endowed with a right action of a group $\Gamma$. We say that \textbf{the quotient $T / \Gamma$ exists} in $C$ if the covariant functor
\begin{align*}
C &\rightarrow \mathrm{Sets} \\
U &\mapsto \mathrm{Hom}_{C}(T, U)^{\Gamma}
\end{align*}
is representable by an object of $C$.
\item[$\triangleright$]Let $T$ be an $S$-scheme. An action of a finite group $\Gamma$ on $T$ is \textbf{admissible} if there exists an affine $\Gamma$-invariant morphism $f : T \rightarrow T'$ such that the canonical morphism $\Ow_{T'} \rightarrow f_* \Ow_{T}$ induces an isomorphism from $\Ow_{T'}$ to $(f_* \Ow_{T})^{\Gamma}$.
\end{itemize}

\begin{prop}[\cite{SGA1}, V.1.3] \label{chap2sga1.1.3}Let $T$ be an $S$-scheme endowed with an admissible right action of a finite group $\Gamma$. If $f : T \rightarrow T'$ is an affine $\Gamma$-invariant morphism such that the canonical morphism $\Ow_{T'} \rightarrow f_* \Ow_{T}$ induces an isomorphism from $\Ow_{T'}$ to $(f_* \Ow_{T})^{\Gamma}$, then the quotient $T/ \Gamma$ exists and is isomorphic to $T'$.
\end{prop}

\begin{prop}[\cite{SGA1}, V.1.8] \label{chap2sga1.1.8}Let $T$ be an $S$-scheme endowed with a right action of a finite group $\Gamma$. Then, the action of $\Gamma$ on $T$ is admissible if and only if $T$ is covered by $\Gamma$-invariant affine open subsets.
\end{prop}

\begin{prop}[\cite{SGA1}, V.1.9] \label{chap2sga1.1.9}Let $T$ be an $S$-scheme endowed with an admissible right action of a finite group $\Gamma$, and let $S'$ be a flat $S$-scheme. Then, the action of $\Gamma$ on the $S'$-scheme $T \times_S S'$ is admissible, and the canonical morphism
$$
(T \times_S S')/ \Gamma \rightarrow (T / \Gamma) \times_S S'
$$
is an isomorphism.
\end{prop}

Let $X$ be an $S$-scheme and let $d \geq 0$ be an integer. The group $\mathfrak{S}_d$ of permutations of $\llbracket1,d\rrbracket$ acts on the right on the $S$-scheme $X^{\times_S d} = X \times_S \dots \times_S X$ by the formula
$$
(x_i)_{i \in \llbracket1,d\rrbracket} \cdot \sigma = (x_{\sigma(i)})_{i \in \llbracket1,d\rrbracket}.
$$

\begin{prop}\label{chap2symdrepr}If $X$ is Zariski-locally quasi-projective over $S$, then the right action of $\mathfrak{S}_d$ on $X^{\times_S d}$ is admissible. In particular, the quotient $\Sy_S^d(X) = X^{\times_S d} / \mathfrak{S}_d $  exists in the category of $S$-schemes.
\end{prop}

Since $X$ is Zariski-locally quasi-projective over $S$, any finite set of points in $X$ with the same image in $S$ is contained in an affine open subset of $X$. Thus $X^{\times_S d}$ is covered by open subsets of the form $U^{\times_S d}$ where $U$ is an affine open subset of $X$ whose image in $S$ is contained in an affine open subset of $S$. These particular open subsets are affine and $\mathfrak{S}_d$-invariant, so that the action of $\mathfrak{S}_d$ on $X^{\times_S d}$ is admissible by Proposition \ref{chap2sga1.1.8}.

\begin{rema}\label{chap2rema15}If $X = \Spec(B)$ and $S = \Spec(A)$ are affine, then for any $S$-scheme $T$ we have 
\begin{align*}
\mathrm{Hom}_{\Sch_{/S}}(  X^{\times_S d} , T)^{\mathfrak{S}_d} &= \mathrm{Hom}_{\mathrm{Alg}_A}( \Gamma(T,\Ow_T), B^{\otimes_A d})^{\mathfrak{S}_d} \\
&= \mathrm{Hom}_{\mathrm{Alg}_A}( \Gamma(T,\Ow_T),\mathrm{TS}_A^d(B)),
\end{align*}
cf. \ref{chap23.10}. Thus $\Sym^d_{S}(X)$ is representable by the $S$-scheme $\Spec(\mathrm{TS}_A^d(B))$.
\end{rema}

\begin{prop}\label{chap2symdrepr2}If $X$ is flat and Zariski-locally quasi-projective over $S$, then $\Sy_S^d(X)$ is flat over $S$. Moreover, for any $S$-scheme $S'$, the canonical morphism
$$
\Sy_{S'}^d(X \times_{S} S') \rightarrow \Sy_S^d(X) \times_S S'
$$
is an isomorphism.
\end{prop}

This follows from Remark \ref{chap2rema15} and from Proposition \ref{chap2flat}.

\begin{prop}[\cite{SGA1}, IX.5.8]\label{chap2torsdescent2} Let $G$ be a finite abelian group, let $P$ be a $G$-torsor over an $S$-scheme $X$ in $S_{\et}$. Assume that $P$ and $X$ are endowed with right actions from a finite group $\Gamma$ such that the morphism $P \rightarrow X$ is $\Gamma$-equivariant, and that the following properties hold:
\begin{itemize}
\item[$(a)$] The right $\Gamma$-action on $P$ commutes with the left $G$-action.
\item[$(b)$] The right $\Gamma$-action on $X$ is admissible (cf. \ref{chap2defquot}), and the quotient morphism $X \rightarrow X/ \Gamma$ is finite.
\item[$(c)$] For any geometric point $\bar{x}$ of $X$, the action of the stabilizer $\Gamma_{\bar{x}}$ of $\bar{x}$ in $\Gamma$ on the fiber $P_{\bar{x}}$ of $P$ at $\bar{x}$ is trivial.
\end{itemize}
Then the action of $\Gamma$ on $P$ is admissible, and $P/\Gamma$ is a $G$-torsor over $X/\Gamma$ in $S_{\et}$.
\end{prop}

\subsection{\label{chap23.12}} Let $S$ be a scheme, let $X$ be an $S$-scheme and let $d \geq 1$ be an integer. Let $G$ be a finite abelian group, and let $P \rightarrow X$ be a $G$-torsor over $X$ in $S_{\et}$. By \ref{chap2torsrepre}, the sheaf $P$ is representable by a finite \'etale $X$-scheme.

For each $i \in \llbracket1,d\rrbracket$ let $p_i : X^{\times_S d} \rightarrow X$ be the projection on $i$-th factor, and let us consider the $G$-torsor
$$
p_1^{-1} P \otimes \dots \otimes p_d^{-1} P = G_d \setminus P^{\times_S d}
$$
over $X^{\times_S d}$, where $G_d \subseteq G^{ d}$ is the kernel of the multiplication morphism $G^{ d} \rightarrow G$. By \ref{chap2torsrepre}, the object $G_d \setminus P^{\times_S d}$ of $S_{\et}$ is representable by an $S$-scheme which is finite \'etale over $X^{\times_S d}$. The group $\mathfrak{S}_d$ acts on the right on $G_d \setminus P^{\times_S d}$ by the formula
$$
(p_i)_{i \in \llbracket1,d\rrbracket} \cdot \sigma = (p_{\sigma(i)})_{i \in \llbracket1,d\rrbracket}.
$$
This action of $\mathfrak{S}_d$ commutes with the left action of $G$ on $G_d \setminus P^{\times_S d}$.

\begin{prop}\label{chap2descent} If $X$ is Zariski-locally quasi-projective on $S$, then the right action of $\mathfrak{S}_d$ on $G_d \setminus P^{\times_S d}$ is admissible (cf. \ref{chap2defquot}), so that the quotient $P^{[d]}$ of $G_d \setminus P^{\times_S d}$ by $\mathfrak{S}_d$ exists in $\Sch_{/S}$. Moreover, the canonical morphism $P^{[d]} \rightarrow \Sym^d_{S}(X)$ is a $G$-torsor, and the morphism
$$
p_1^{-1} P \otimes \dots \otimes p_d^{-1} P \rightarrow r^{-1} P^{[d]}
$$
where $r : X^{\times_S d} \rightarrow \Sym^d_{S}(X)$ is the canonical projection, is an isomorphism of $G$-torsors over $X^{\times_S d}$.
\end{prop}

By \ref{chap2symdrepr} and \ref{chap2torsdescent2}, it is sufficient to show that if $\bar{x} = (\bar{x}_i)_{i=1}^d$ is a geometric point of $X^{\times_S d}$, then the stabilizer of $\bar{x}$ in $\mathfrak{S}_d$ acts trivially on $(G_d \setminus P^{\times_S d} )_{\bar{x}}$. 
Assume that the finite set $\{ \bar{x}_i \ | \ i \in \llbracket1,d] \} $ has exactly $r$ distinct elements $\bar{y}_1,\dots,\bar{y}_r$, where $\bar{y}_j$ appears with multiplicity $d_j$. Then the stabilizer of $\bar{x}$ in $\mathfrak{S}_d$ is isomorphic to the subgroup $\prod_{j=1}^r \mathfrak{S}_{d_j}$ of $\mathfrak{S}_d$. For each $j \in \llbracket1,r\rrbracket$, the $G$-torsor $ P_{\bar{y}_j}$ is trivial, and if $e$ is a section of this torsor then $(e)_{i=1}^{d_j}$ is a section of $G_{d_j} \setminus P_{\bar{y}_j}^{ d_j}$ which is $\mathfrak{S}_{d_j}$-invariant. The action of $\mathfrak{S}_{d_j}$ on $G_{d_j} \setminus P_{\bar{y}_j}^{ d_j}$ is therefore trivial, so that the action of $\prod_{j=1}^r \mathfrak{S}_{d_j}$ on the $G$-torsor
$$
(G_d \setminus P^{\times_S d})_{\bar{x}} = G_r \setminus \left( \prod_{j=1}^r G_{d_j} \setminus P_{\bar{y}_j}^{ d_j} \right)
$$
is trivial as well.

\begin{prop}\label{chap2descent2} If $X$ is flat and Zariski-locally quasi-projective on $S$, then for any $S$-scheme $S'$, the canonical morphism
$$
(P \times_S S')^{[d]} \rightarrow P^{[d]} \times_S S'
$$
is an isomorphism.
\end{prop}

By Proposition \ref{chap2symdrepr2}, the canonical morphism
$$
\Sy_{S'}^d(X \times_{S} S') \rightarrow \Sy_S^d(X) \times_S S'
$$
is an isomorphism. Thus the second morphism in the composition
$$
(P \times_S S')^{[d]} \rightarrow (P^{[d]} \times_S S') \times_{\Sy_S^d(X) \times_S S'} \Sy_{S'}^d(X \times_{S} S')   \rightarrow P^{[d]} \times_S S'
$$
is an isomorphism, while the first morphism is a morphism of $G$-torsors, hence an isomorphism.

\section{Geometric local class field theory \label{chap2glcftpart}}

Let $k$ be a perfect field, and let $L$ be a complete discretely valued extension of $k$ with residue field $k$. We denote by $\Ow_L$ its ring of integers, and by $\mathfrak{m}_L$ the maximal ideal of $\Ow_L$.

\subsection{\label{chap23.2}}Let us consider the functor
\begin{align*}
\mathbb{O}_L : \Alg_{k} &\rightarrow \Alg_{\Ow_L} \\
 A &\mapsto \lim_n A \otimes_k \Ow_L/\mathfrak{m}_L^n,
\end{align*}
with values in the category of $\Ow_L$-algebras.

\begin{prop}\label{chap2repr} The functor $\mathbb{O}_L$ is representable by a $k$-scheme.
\end{prop}

Indeed, if $\pi$ is a uniformizer of $L$, then we have an isomorphism $k((t)) \rightarrow L$ which sends $t$ to $\pi$, so that the functor $\mathbb{O}_L$ is isomorphic to the functor $A \mapsto A[[t]]$, which is representable by an affine space over $k$ of countable dimension.

\begin{cor} The functor $\mathbb{L} = \mathbb{O}_L \otimes_{\Ow_L} L $ is representable by an ind-$k$-scheme.
\end{cor}

We can assume that $L$ is the field of Laurent series $k((t))$. In this case, we have
$$
\mathbb{L}(A) = A((t)) = \mathrm{colim}_n \ t^{-n} A[[t]]
$$
for any $k$-algebra $A$, and for each integer $n$ the functor $A \mapsto t^{-n} A[[t]]$ is representable by a $k$-scheme, cf. \ref{chap2repr}.

\begin{prop}\label{chap2repr2} Let $G$ (resp. $H$) be the functor from $ \Alg_{k}$ to the category of groups which associates to a $k$-algebra $A$ the subgroup $G(A)$ of $A((t))^{\times}$ consisting of Laurent series of the form $1 + \sum_{r>0} a_r t^{-r}$ where $a_r$ is a nilpotent element of $A$ for each $r > 0$ and vanishes for $r$ large enough (resp. of Laurent series of the form $1 + \sum_{r>0} a_r t^{r}$ where $a_r$ belongs to $A$ for each $r > 0$). Let $\underline{\mathbb{Z}}$ be the functor which sends a $k$-algebra $A$ to the group of locally constant functions $\Spec(A) \rightarrow \mathbb{Z}$. Then for any uniformizer $\pi$ of $L$, the morphism
\begin{align*}
 \mathbb{G}_{m,k} \times \underline{\mathbb{Z}} \times G  \times H &\rightarrow \mathbb{L}^{\times},\\
(a,n,g,h) &\mapsto a \pi^n  g(\pi) h(\pi),
\end{align*}
is an isomorphism of group-valued functors.
\end{prop}

Let $A$ be a $k$-algebra. By (\cite{CC}, 0.8), every invertible element $u$ of $A((t))$ uniquely factors as $u = t^n f(t) h(t)$ where $f(t)$ and $h(t)$ are elements of $A[[t]]^{\times}$ and $G(A)$ respectively, and $n : \Spec(A) \rightarrow \mathbb{Z}$ is a locally constant function. Moreover, there is a unique factorisation $f(t) = a g(t)$ where $a$ and $g(t)$ belong to $A^{\times}$ and $H(A)$ respectively, hence the result.

\begin{cor}\label{chap2repr3} The functor $\mathbb{L}^{\times}$ is representable by an ind-$k$-scheme. Moreover, its restriction to the category of reduced $k$-algebras is representable by a reduced $k$-scheme.
\end{cor}

The groups $\underline{\mathbb{Z}}$ and $H$ from Proposition \ref{chap2repr2} are representable by reduced $k$-schemes, and so is $\mathbb{G}_{m,k}$. Moreover, the group $G$ from \ref{chap2repr2} is the filtered colimit of the functor $n \mapsto G_n$, where $G_n$ is the functor which associates to a $k$-algebra $A$ the subset $G_n(A)$ of $A((t))^{\times}$ consisting of Laurent series of the form $1 + \sum_{r=1}^n a_r t^{-r}$ where $a_r^n = 0$ for each $r \in  \llbracket1,n\rrbracket$. For each $n$, the functor $G_n$ is representable by an affine $k$-scheme. Thus $G$ is representable by an ind-$k$-scheme, and so is $\mathbb{L}^{\times}$ by \ref{chap2repr2}. The last assertion of Corollary \ref{chap2repr3} follows from the fact that $G(A)$ is the trivial group for any reduced $k$-algebra $A$.

\begin{cor}\label{chap2repr4} Let $d \geq 0$ be an integer, and let $\mathbb{U}_L^{(d)}$ be the subfunctor $1 + \mathfrak{m}_L^d \mathbb{O}_L$ (resp. $\mathbb{O}_L^{\times}$) of $\mathbb{L}^{\times}$ for $d \geq 1$ (resp. for $d=0$). Then the functor
\begin{align*}
\mathbb{L}^{\times} / \mathbb{U}_L^{(d)} : \Alg_{k} &\rightarrow \mathrm{Sets} \\
 A  &\mapsto \mathbb{L}^{\times}(A) / \mathbb{U}_L^{(d)}(A),
\end{align*}
is representable by an ind-$k$-scheme. Moreover, its restriction to the category of reduced $k$-algebras is representable by a reduced $k$-scheme.
\end{cor}

According to Proposition \ref{chap2repr2}, it is sufficent to show that $(\mathbb{G}_{m,k} \times H)/ \mathbb{U}_{k((t))}^{(d)}$ is representable by a reduced $k$-scheme. The case $d = 0$ is clear, while for $d \geq 1$, we have for any $k$-algebra $A$ a bijection
\begin{align*}
A^{\times} \times A^{\llbracket1,d-1\rrbracket} &\rightarrow (\mathbb{G}_{m,k} \times H)(A)/ \mathbb{U}_{k((t))}^{(d)}(A) \\
 (a_i)_{0 \leq i \leq d-1}  &\mapsto \sum_{i=0}^{d-1} a_i t^i,
\end{align*}
hence the result.

\subsection{\label{chap23.3}}From now on, we consider $\Spec(L)$, $\mathbb{L}^{\times}$ and $\mathbb{L}^{\times} / \mathbb{U}_L^{(d)}$ for each integer $d \geq 0$ as objects of the topos $\Spec(k)_{\et}$. 
Let $\pi$ be an uniformizer of $L$. We denote by $\Pi$ the element of $\mathbb{L}(k)$ corresponding to $\pi$ via the canonical identification $L \simeq \mathbb{L}(k)$ . Thus the functor $\mathbb{L}^{\times}$ is given by
$$
\mathbb{L}^{\times} : A \in \Alg_{k}  \mapsto A((\Pi))^{\times}.
$$
In particular, the Laurent series $(\Pi-\pi)^{-1} \Pi = - \sum_{n \geq 1} \pi^{-n} \Pi^n$ defines an $L$-point of $\mathbb{L}^{\times}$. We denote by $\varphi : \Spec(L) \rightarrow \mathbb{L}^{\times}$ the corresponding morphism. We follow here Contou-Carr\`ere's convention; in \cite{TS}, the morphism $\varphi$ corresponds to the point $(\Pi-\pi) \Pi^{-1}$ instead. This is harmless since the inversion is an automorphism of the abelian group $\mathbb{L}^{\times}$.

\begin{teo}[\cite{TS}, Th. A (1)]\label{chap2lcft} Let $G$ be a finite abelian group. The functor
\begin{align*}
\mathrm{Tors}^{\otimes}(\mathbb{L}^{\times},G) &\rightarrow  \mathrm{Tors}(\Spec(L),G)\\
P &\rightarrow \varphi^{-1} P
\end{align*}
is an equivalence of categories (cf. \ref{chap2defitorseur}, \ref{chap2chartors}).
\end{teo}

In the case where $k$ is algebraically closed, Serre constructed in \cite{JPS} an equivalence
$$
\mathrm{Tors}(\Spec(L),G) \rightarrow \mathrm{Tors}^{\otimes}(\mathbb{L}^{\times},G).
$$
More precisely, Serre considers \'etale isogenies over $\mathbb{L}^{\times}$ and the link with $\mathrm{Tors}^{\otimes}(\mathbb{L}^{\times},G)$ is provided by \ref{chap23.1}.
In \cite{TS}, Suzuki shows that the functor from Theorem \ref{chap2lcft} is a quasi-inverse to Serre's functor when $k$ is algebraically closed, and extends the result to arbitrary perfect residue fields. In particular, the equivalence from Theorem \ref{chap2lcft} is canonical, even though its definition depends on the choice of $\pi$. Suzuki's proof of Theorem \ref{chap2lcft} relies on the Albanese property of the morphism $\varphi$, previously established by Contou-Carr\`ere.

Let $L^{\mathrm{sep}}$ be a separable closure of $L$, and let $G_L$ be the Galois group of $L^{\mathrm{sep}}$ over $L$, so that the small \'etale topos of $\Spec(L)$ is isomorphic to the topos of sets with continuous left $G_L$-action. By \ref{chap22.13}, the category of $G$-torsors over $\Spec(L)$ in $\Spec(k)_{\et}$ is isomorphic to the category of $G$-torsors in the small \'etale topos $\Spec(L)_{\eet}$. Correspondingly, for each finite abelian group $G$, the group of isomorphism classes of the category $\mathrm{Tors}(\Spec(L),G)$ is isomorphic to the group of continuous homomorphisms from $G_L$ to $G$.

We denote by $(G_L^j)_{j \geq -1}$ the ramification filtration of $G_L$ (\cite{JPS2}, IV.3), so that $G_L^{-1} = G_L$ and $G_L^{0}$ is the inertia subgroup of $G_L$, while $G_L^{0+} = \cup_{j >0} G_L^j$ is the wild inertia subgroup of $G_L$.

\begin{defi}\label{chap2ramif} Let $G$ be a finite abelian group and let $d \geq 0$ be a rational number. A $G$-torsor over $\Spec(L)$ (in $\Spec(k)_{\et}$), corresponding to a continuous homomorphism $\rho : G_L \rightarrow G$, is said to have \textbf{ramification bounded by $d$} if $\rho(G_L^d) = \{ 1 \}$. A $G$-torsor over $\Spec(L)$ with ramification bounded by $0$ (resp. $1$) is said to be unramified (resp. tamely ramified).
 \end{defi}

\begin{rema} If $P \rightarrow \Spec(L)$ is a $G$-torsor in $\Spec(k)_{\et}$, then we have a finite decomposition
$$
P = \coprod_{i \in I} \Spec(L_i),
$$
where each $L_i$ is a finite separable extension of $L$, and are pairwise isomorphic. The $G$-torsor $P$ has ramification bounded by $d$ if and only if for each $i$ (or, equivalently, for some $i$) the extension $L_i/L$ has ramification bounded by $d$, in the sense $G_L^d$ acts trivially on the finite set $\mathrm{Hom}_L(L_i,L^{\mathrm{sep}})$.
\end{rema}

\begin{prop}\label{chap2ramif7} Let $G$ be a finite abelian group, let $d \geq 0$ be an integer, and let $P$ be a multiplicative $G$-torsor $P$ over $\mathbb{L}^{\times}$ (cf. \ref{chap2chartors}). Assume that $k$ is algebraically closed. Then $\varphi^{-1} P$ has ramification bounded by $d$ (cf. \ref{chap2ramif}) if and only if $P$ is the pullback of a multiplicative $G$-torsor over $\mathbb{L}^{\times}/\mathbb{U}_L^{(d)}$ (cf. \ref{chap2repr4}).
\end{prop}

This follows from (\cite{JPS}, 3.2 Th. 1) and from the compatibility of $\varphi^{-1}$ with Serre's construction (\cite{TS}, Th. A (2)).

\subsection{\label{chap23.5}} Let $\pi$ and $\varphi$ be as in \ref{chap23.3}. Let $K$ be a closed sub-extension of $k$ in $L$, such that $K \rightarrow L$ is a finite extension of degree $d$. Since $L$ is a finite free $K$-algebra of rank $d$, we have a canonical morphism of $K$-schemes
$$
\psi : \Spec(K) \rightarrow \Sym_K^d(\Spec(L))
$$
by \ref{chap2morphalg}.

\begin{prop}\label{chap2pointeis} The composition
$$
\Spec(K) \xrightarrow[]{\psi} \Sym_K^d(\Spec(L)) \rightarrow \Sym_k^d(\Spec(L)) \xrightarrow[]{\Sym_k^d(\varphi)} \Sym_k^d(\mathbb{L^{\times}}) \rightarrow \mathbb{L^{\times}},
$$
where the last morphism is given by the multiplication, corresponds to the $K$-point $P_{\pi}(\Pi)^{-1} \Pi^d$ of $\mathbb{L}^{\times}$, where the polynomial $P_{\pi}$ is the characteristic polynomial of the $K$-linear endomorphism $x \mapsto \pi x$ of $L$.
\end{prop}

We first describe the morphism $\psi$. The scheme $\Sym_K^d(\Spec(L))$ is the spectrum of the $k$-algebra $\mathrm{TS}_K^d(L)$ of symmetric tensors of degree $d$ in $L$, cf. \ref{chap2symdrepr}. The elements $e_i = \pi^{i-1}$ for $i=1,\dots,d$ form a $K$-basis of $L$, so that we have a decomposition
$$
\mathrm{TS}_K^d(L) = \bigoplus_{\substack{\alpha : \llbracket 1,d\rrbracket \rightarrow \mathbb{N} \\ \sum_{i} \alpha(i) = d }} K e_{\alpha},
$$
where we have set (cf. \ref{chap23.10})
$$
e_{\alpha} = \sum_{\substack{\beta : \llbracket 1,d\rrbracket \rightarrow \llbracket 1,d\rrbracket \\ \forall i, |\beta^{-1}(\{i\})| = \alpha(i) }} e_{\beta(1)} \otimes \dots \otimes e_{\beta(d)}.
$$
Let us write the norm polynomial as
$$
N_{L/K}\left(\sum_{i=1}^d x_i e_i \right) = \sum_{\substack{\alpha : \llbracket 1,d\rrbracket \rightarrow \mathbb{N} \\ \sum_{i} \alpha(i) = d }} f_{\alpha} x^{\alpha},
$$
where $x^{\alpha} = x_1^{\alpha(1)} \dots x_d^{\alpha(d)}$, and the $f_{\alpha}$'s are uniquely determined elements of $K$. The morphism $\mathrm{TS}_K^d(L) \rightarrow K$ corresponding to $\psi$ is the unique $K$-linear homomorphism which sends $e_{\alpha}$ to $f_{\alpha}$ (cf. \ref{chap2flat2} and its proof).

Next we describe the composition 
$$
\Sym_K^d(\Spec(L)) \rightarrow \Sym_k^d(\Spec(L)) \xrightarrow[]{\Sym_k^d(\varphi)} \Sym_k^d(\mathbb{L^{\times}}) \rightarrow \mathbb{L^{\times}}.
$$
Its precomposition with the projection $\Spec(L)^{\times_K d} \rightarrow \Sym_K^d(\Spec(L))$ corresponds to the element of $L^{\otimes_K d}((\Pi))^{\times}$ given by the formula
$$
\prod_{i=1}^d \left( (\Pi - 1^{\otimes (i-1)} \otimes \pi \otimes 1^{\otimes (d-i)})^{-1} \Pi \right) = P(\Pi)^{-1} \Pi^d,
$$
where the polynomial $P(\Pi)$ can be computed as follows:
\begin{align*}
P(\Pi) &= \prod_{i=1}^d (\Pi - 1^{\otimes (i-1)} \otimes \pi \otimes 1^{\otimes (d-i)}) \\
&= \sum_{r=0}^d (-1)^r \Pi^{d-r} \sum_{\substack{(i_1,\dots,i_d)\in \{0,1 \}^d \\ |\{ s | i_s = 1 \}| = r }} \pi^{i_1} \otimes \dots \otimes \pi^{i_d} \\
&= \sum_{r=0}^d (-1)^r e_{\alpha_r}\Pi^{d-r},
\end{align*}
where $\alpha_r :\llbracket 1,d\rrbracket \rightarrow \mathbb{N} $ is the map which sends $1$ and $2$ to $d-r$ and $r$ respectively, and any $i > 2$ to $0$.
The image of $P(\Pi)$ by $\psi$ in $K[\Pi]$ is the polynomial
$$
\sum_{r=0}^d (-1)^r f_{\alpha_r}\Pi^{d-r} = N_{L[\Pi]/K[\Pi]} \left( \Pi e_1 - e_2 \right).
$$
Since $e_1 =1$ and $e_2 = \pi$, we obtain \ref{chap2pointeis}.

\begin{prop}\label{chap2ramif2} Let $G$ be a finite abelian group, and let $Q$ be a $G$-torsor over $\Spec(L)$ (in $\Spec(k)_{\et}$) of ramification bounded by $d$ (cf. \ref{chap2ramif}). Then $\psi^{-1} Q^{[d]}$ (cf. \ref{chap2descent}) is tamely ramified on $\Spec(K)$.
\end{prop}

Let $K'$ be the maximal unramified extension of $K$ in a separable closure of $K$. The formation of $\Sym_K^d(\Spec(L))$ is compatible with any base change by Proposition \ref{chap2sga1.1.9} or by Proposition \ref{chap2symdrepr2}, and so is the formation of $\varphi$. Moreover, a $G$-torsor over $\Spec(K)$ is tamely ramified if and only if its restriction to $\Spec(K')$ is tamely ramified. By replacing $K$ and $L$ by $K'$ and the components of $K' \otimes_K L$ respectively, we can assume that the residue field $k$ is algebraically closed.

Let $P$ be the multiplicative $G$-torsor on $\mathbb{L}^{\times}$  (cf. \ref{chap2chartors}) associated to $Q$ (cf. \ref{chap2lcft}), so that $Q$ is isomorphic to $\varphi^{-1} P$. Then $\psi^{-1} Q^{[d]}$ is isomorphic to the pullback of $P$ along the composition 
$$
\Spec(K) \xrightarrow[]{\psi} \Sym_K^d(\Spec(L)) \rightarrow \Sym_k^d(\Spec(L)) \xrightarrow[]{\Sym_k^d(\varphi)} \Sym_k^d(\mathbb{L^{\times}}) \rightarrow \mathbb{L^{\times}}
$$
considered in \ref{chap2pointeis}. 
By \ref{chap2pointeis}, this composition corresponds to the $K$-point of $\mathbb{L}^{\times}$ given by $P_{\pi}(\Pi)^{-1} \Pi^d$, where $P_{\pi}$ is the characteristic polynomial of $\pi$ acting $K$-linearly by multiplication on $L$. Let us consider the morphism of pointed sets
\begin{align*}
\rho : \mathbb{L}^{\times}(K) &\rightarrow H^1(\Spec(K)_{\et},G) \\
\nu &\rightarrow \nu^{-1} P
\end{align*}
where an element $\nu$ of $\mathbb{L}^{\times}(K)$ is identified to a morphism $\Spec(K) \rightarrow \mathbb{L}^{\times} $. If $\nu_1$ and $\nu_2$ are elements of $\mathbb{L}^{\times}(K)$, then using the isomorphism $\theta : p_1^{-1} P \otimes p_2^{-1} P \rightarrow m^{-1} P$ from \ref{chap2chartors}, we obtain isomorphisms
$$
(\nu_1 \nu_2)^{-1} P  \leftarrow (\nu_1 \times \nu_2)^{-1} m^{-1} P \xleftarrow[]{(\nu_1 \times \nu_2)^{-1}\theta} (\nu_1 \times \nu_2)^{-1}(p_1^{-1} P \otimes p_2^{-1} P )\leftarrow \nu_1^{-1} P \otimes \nu_2^{-1} P.
$$
Thus $\rho$ is an homomorphism of abelian groups.

We have to prove that $\rho(\nu)$ is the isomorphism class of a tamely ramified $G$-torsor over $\Spec(K)$, where $\nu = P_{\pi}(\Pi)^{-1} \Pi^d $. Since $P_{\pi}$ is an Eisenstein polynomial, it can be written as $P_{\pi}(\Pi) = \Pi^d + c R(\Pi)$, where $c = P_{\pi}(0)$ is a uniformizer of $K$, and $R$ is a polynomial of degree $<d$ with coefficients in $\Ow_K$, such that $R(0)=1$. Thus we can write
$$
\nu = c^{-1} \nu_1 \nu_2,
$$
where $\nu_1 = R(\Pi)^{-1} \Pi^d$ and $\nu_2 = (1 + c^{-1} \Pi^d R(\Pi)^{-1})^{-1}$, so that $\rho(\nu) = \rho(c)^{-1} \rho(\nu_1) \rho(\nu_2)$.

Since $Q$ has ramification bounded by $d$ (cf. \ref{chap2ramif}), the restriction of $\rho$ to $\mathbb{U}_L^{(d)}(K)$ is trivial (cf. \ref{chap2ramif7}). In particular, $\rho(\nu_2)$ is trivial since $\nu_2$ belongs to $\mathbb{U}_L^{(d)}(K)$.

The element $\nu_1$ belongs to $ \mathbb{L}^{\times}(\Ow_K)$, so that the morphism $\nu_1 : \Spec(K) \rightarrow \mathbb{L}^{\times} $ factors through $\Spec(\Ow_K)$. This implies that $\rho(\nu_1)$ is the isomorphism class of an unramified $G$-torsor over $\Spec(K)$. It remains to prove that $\rho(c)$ is the isomorphism class of a tamely ramified $G$-torsor over $\Spec(K)$. Since $c$ belongs to $K^{\times} = \mathbb{G}_{m,k}(K) \subseteq \mathbb{L}^{\times}(K)$, this is a consequence of the following lemma:

\begin{lem} Let $T$ be a multiplicative $G$-torsor over the $k$-group scheme $\mathbb{G}_{m,k}$ (cf. \ref{chap2chartors}). Then $T$ is tamely ramified at $0$ and $\infty$. 
\end{lem}

Let $G_k$ be the constant $k$-group scheme associated to $k$. By \ref{chap2ext}, there is a structure of $k$-group scheme on $T$ and an exact sequence
\begin{align}\label{chap2ext10}
1 \rightarrow G_k \rightarrow T \rightarrow \mathbb{G}_{m,k} \rightarrow 1
\end{align}
in $\Spec(k)_{\et}$, such that the structure of $G$-torsor on $T$ is given by the action of its subgroup $G$ by translations. Since the $\fppf$ topology is finer than the \'etale topology on $\Sch_{/k}$, the sequence \ref{chap2ext10} remains exact in the topos $\Spec(k)_{\Fppf}$. In particular, we obtain a class in the group $\mathrm{Ext}^1_{\Fppf}(\mathbb{G}_{m,k},G_k )$ of extensions of $\mathbb{G}_{m,k}$ by $G_k$ in $\Spec(k)_{\Fppf}$.

Let $n = |G|$. In the topos $\Spec(k)_{\Fppf}$ we have an exact sequence
\begin{align}\label{chap2ext11}
1 \rightarrow \mu_{n,k} \rightarrow \mathbb{G}_{m,k} \xrightarrow[]{n} \mathbb{G}_{m,k} \rightarrow 1,
\end{align}
where $\mu_{n,k}$ is the $k$-group scheme of $n$-th roots of unity. By applying the functor $\mathrm{Hom}(\cdot, G_k)$, we obtain an exact sequence
$$
\mathrm{Hom}( \mu_{n,k} , G_k) \xrightarrow[]{\delta} \mathrm{Ext}^1_{\fppf}(\mathbb{G}_{m,k},G_k )  \xrightarrow[]{n} \mathrm{Ext}^1_{\fppf}(\mathbb{G}_{m,k},G_k ).
$$
Since $n = |G|$, the group $\mathrm{Ext}^1_{\Fppf}(\mathbb{G}_{m,k},G_k )$ is annihilated by $n$, so that the homomorphism $\delta$ above is surjective. Thus the exact sequence (\ref{chap2ext10}) in $\Spec(k)_{\Fppf}$ is the pushout of (\ref{chap2ext11}) along an homomorphism $\mu_{n,k} \rightarrow G_k$. Let $n'$ be the largest divisor of $n$ which is invertible in $k$. Then the largest \'etale quotient of $\mu_{n,k}$ is the epimorphism $\mu_{n,k} \rightarrow \mu_{n',k}$ given by $x \mapsto x^{\frac{n}{n'}}$. In particular, the homomorphism $\mu_{n,k} \rightarrow G_k$ factors through $\mu_{n',k}$, so that (\ref{chap2ext10}) is the pushout of the extension 
$$
1 \rightarrow \mu_{n',k} \rightarrow \mathbb{G}_{m,k}  \xrightarrow[]{n'}  \mathbb{G}_{m,k} \rightarrow  1
$$
along an homomorphism $\mu_{n',k} \rightarrow G_k$. Since the morphism $\mathbb{G}_{m,k}  \xrightarrow[]{n'}  \mathbb{G}_{m,k}$ is tamely ramified above $0$ and $\infty$, so is the morphism $T \rightarrow \mathbb{G}_{m,k}$.

\section{Rigidified Picard schemes of relative curves \label{chap2rigpicpart}}

\subsection{\label{chap21.0}} Let $f : X \rightarrow S$ be a smooth morphism of schemes of relative dimension $1$, with connected geometric fibers of genus $g$, which is Zariski-locally projective over $S$. 

\begin{prop}\label{chap2coh} The canonical homomorphism $\Ow_S \rightarrow f_* \Ow_X$ is an isomorphism.
\end{prop}

If $S$ is locally noetherian, then $\Ow_X$ is cohomologically flat over $S$ in dimension $0$ by (\cite{EGA3}, 7.8.6). This means that for any quasi-coherent $\Ow_S$-module $\M$, the canonical homomorphism $f_* f^* \Ow_X \otimes_{\Ow_S} \M \rightarrow f_* f^* \M$ is an isomorphism. This implies that the formation of $f_* \Ow_X$ commutes with arbitrary base change: if $f' : X \times_S S' \rightarrow S'$ is the base change of $f$ by a morphism of schemes $S' \rightarrow S$, then the canonical morphism $f_* \Ow_{X} \otimes_{\Ow_S} \Ow_{S'} \rightarrow f'_* \Ow_{ X \times_S S'}$ is an isomorphism, cf. (\cite{EGA3}, 7.7.5.3). By applying this result to the inclusion $\Spec(\ka(s)) \rightarrow S$ of a point $s$ of $S$, we obtain that $f_*(\Ow_X)_s \otimes_{\Ow_{S,s}} \ka(s)$ is isomorphic to  $H^0(X_s,\Ow_{X_s}) = \ka(s)$. Since $f_*(\Ow_X)$ is a coherent $\Ow_S$-module, Nakayama's lemma yield that the canonical morphism $\Ow_S \rightarrow f_*(\Ow_X)$ is an epimorphism. It is also injective since $f$ is faithfully flat, hence the result.

In general one can assume that $S$ is affine and that $X$ is projective over $S$, in which case there is a noetherian scheme $S_0$, a morphism $S \rightarrow S_0$ and a smooth projective $S_0$-scheme $X_0$ with geometrically connected fibers such that $X$ is isomorphic to the $S$-scheme $X_0 \times_{S_0} S$, cf. (\cite{EGA42}, 8.9.1, 8.10.5(xiii), 17.7.9). We have already seen that in this case the canonical homomorphism $\Ow_{S_0} \rightarrow f_* \Ow_{X_0}$ is an isomorphism, and that the formation of $f_* \Ow_{X_0}$ commutes with arbitrary base change. In particular, both morphisms in the sequence
$$
\Ow_{S} \rightarrow f_* \Ow_{X_0} \otimes_{\Ow_{S_0}} \Ow_{S} \rightarrow f_* \Ow_{ X}
$$ 
are isomorphisms.

\begin{prop}\label{chap2coh2} Let $d \geq 2g - 1$ be an integer, and let $\Lc$ be an invertible $\Ow_X$-module with degree $d$ on each fiber of $f$. Then, the $\Ow_S$-module $f_* \Lc$ is locally free of rank $d-g+1$, the higher direct images $R^j(f_* \Lc)$ vanish for $j >0$, and the formation of $f_* \Lc$ commutes with arbitrary base change: if $f' : X' \rightarrow S'$ is the base change of $f$ by a morphism $S' \rightarrow S$, then the canonical homomorphism $f_* \Lc \otimes_{\Ow_S} \Ow_{S'} \rightarrow f'_*(\Lc \otimes_{\Ow_X} \Ow_{X'})$ is an isomorphism.
\end{prop}

We first assume that $S$ is locally noetherian. For each point of $s$ of $S$ and for each integer $i$, the Riemann-Roch theorem for smooth projective curves implies that the $k(s)$-vector space $H^i(X_s, \Lc_s)$ is of dimension $d-g+1$ for $i=0$, and vanishes otherwise. This implies that $R^j f_*( \Lc \otimes_{\Ow_X} f^* \mathcal{N})$ vanishes for any integer $j>0$ and any $\Ow_S$-module $\mathcal{N}$ by the proof of (\cite{EGA3}, 7.9.8). Let 
$$0 \rightarrow \mathcal{N} \rightarrow \mathcal{M} \rightarrow \mathcal{P} \rightarrow 0$$
be an exact sequence of $\Ow_S$-modules. Since $f$ is flat and since $\Lc$ is a flat $\Ow_X$-module, the sequence
$$
0 \rightarrow  \Lc \otimes_{\Ow_X}  f^*\mathcal{N} \rightarrow \Lc \otimes_{\Ow_X} f^* \mathcal{M} \rightarrow  \Lc \otimes_{\Ow_X} f^*\mathcal{P} \rightarrow 0
$$
is exact as well. Since $R^1 f_*(\Lc \otimes_{\Ow_X}  f^*\mathcal{N}) $ vanishes, the sequence
$$
0 \rightarrow  f_*(\Lc \otimes_{\Ow_X}  f^*\mathcal{N} ) \rightarrow f_*(\Lc \otimes_{\Ow_X} f^*\mathcal{M}) \rightarrow  f_*( \Lc \otimes_{\Ow_X} f^*\mathcal{P}) \rightarrow 0
$$
is exact. The $\Ow_X$-module $\Lc$ is therefore cohomologically flat over $S$ in dimension $0$, cf. (\cite{EGA3}, 7.8.1). By (\cite{EGA3}, 7.8.4(d)) the $\Ow_S$-module $f_* \Lc$ is locally free, and the formation of $f_* \Lc$ commutes with arbitrary base change. By applying the latter result to the inclusion $\Spec(\ka(s)) \rightarrow S$ of a point $s$ of $S$ and by using that $H^0(X_s, \Lc_s)$ is of dimension $d-g+1$ over $\ka(s)$, we obtain that the locally free $\Ow_X$-module $f_* \Lc$ is of constant rank $d-g+1$.

In general one can assume that $S$ is affine and that $X$ is projective over $S$, in which case there is a noetherian scheme $S_0$, a morphism $S \rightarrow S_0$, a smooth projective $S_0$-scheme $X_0$, and an invertible $\Ow_{X_0}$-module $\Lc_0$ such that $X$ is isomorphic to the $S$-scheme $X_0 \times_{S_0} S$ and $\Lc$ is isomorphic to the pullback of $\Lc_0$ by the canonical projection $X_0 \times_{S_0} S \rightarrow X_0$, cf. (\cite{EGA42}, 8.9.1, 8.10.5(xiii), 17.7.9). We have seen that the $\Ow_{S_0}$-module $f_{0*} \Lc$ is locally free of rank $d-g+1$, and that its formation commutes with arbitrary base change. By performing the base change by the morphism $S \rightarrow S_0$, we obtain that $f_* \Lc$ is a locally free $\Ow_S$-module of rank $d-g+1$ and that the formation of $f_* \Lc$ commutes with arbitrary base change.

\subsection{\label{chap21.0.1}} Let $f : X \rightarrow S$ be as in \ref{chap21.0}. The \textbf{relative Picard functor} of $f$ is the sheaf of abelian groups $\Pic_{S}(X) = R^1 f_{\Fppf,*} \mathbb{G}_m$ in $S_{\Fppf}$. Alternatively, $\Pic_{S}(X)$ is the sheaf of abelian groups on $S$ associated to the presheaf which sends an $S$-scheme $T$ to $\Pic(X \times_S T)$, the abelian group of isomorphism classes of invertible $\Ow_{X \times_S T}$-modules. For any $S$-scheme $S'$, we have $(S_{\Fppf})_{/S'} = S'_{\Fppf}$, and we thus have:
\begin{prop} For any $S$-scheme $S'$, the canonical morphism
$$
\Pic_{S'}(X \times_S S') \rightarrow \Pic_{S}(X) \times_S S'
$$
is an isomorphism in $S'_{\Fppf}$.
\end{prop}

The elements of $\Pic(X \times_S T)$ which are pulled back from an element of $\Pic(T)$ yield trivial classes in $\Pic_S(X)(T)$, since invertible $\Ow_{T}$-modules are locally trivial on $T$ (for the Zariski topology, and thus for the $\fppf$-topology). This yields a sequence
\begin{align}\label{chap2exactseq}
0 \rightarrow \Pic(T) \rightarrow \Pic(X \times_S T) \rightarrow \Pic_S(X)(T) \rightarrow 0,
\end{align}
which is however not necessarily exact. The following is Proposition $4$ from (\cite{BLR}, 8.1), whose assumptions are satisfied by \ref{chap2coh}:

\begin{prop}\label{chap2exactseq2} If $f$ has a section, then the sequence \eqref{chap2exactseq} is exact for any $S$-scheme $T$.
\end{prop}

By a theorem of Grothendieck (\cite{BLR}, 8.2.1) the sheaf $\Pic_S(X)$ is representable by a separated $S$-scheme. By (\cite{BLR}, 9.3.1) the $S$-scheme $\Pic_S(X)$ is smooth of relative dimension $g$, and there is a decomposition 
$$
\Pic_{S}(X) = \coprod_{d \in \Z} \Pic_{S}^d(X),
$$
into open and closed subschemes, where $\Pic_{S}^d(X)$ is the $\fppf$-sheaf associated to the presheaf
\begin{align*}
\Sch_{/S}^{\fp} &\rightarrow \mathrm{Sets}\\
T  &\mapsto \{ \Lc \in \Pic(X \times_S T) | \forall \bar{t} \rightarrow T, \mathrm{deg}_{X_{\bar{t}}}(\Lc_{\bar{t}}) = d \}.
\end{align*}
Here the condition $\mathrm{deg}_{X_{\bar{t}}}(\Lc_{\bar{t}}) = d$ runs over all geometric points $\bar{t} \rightarrow T$ of $T$. 
%

\subsection{\label{chap21.1}}  Let $f : X \rightarrow S$ be as in \ref{chap21.0}, and let $i : Y \hookrightarrow X$ be a closed subscheme of $X$, which is finite locally free over $S$ of degree $N \geq 1$. A \textbf{$Y$-rigidified line bundle on $X$} is a pair $(\Lc,\alpha)$ where $\Lc$ is a locally free $\Ow_X$-module of rank $1$ and $\alpha : \Ow_Y \rightarrow  i^* \Lc$ is an isomorphism of $\Ow_Y$-modules. Two $Y$-rigidified line bundles $(\Lc,\alpha)$ and $(\Lc',\alpha')$ are \textbf{equivalent} if there is an isomorphism $\beta : \Lc \rightarrow \Lc'$ of $\Ow_X$-modules such that $(i^*\beta) \alpha = \alpha'$.
If such an isomorphism $\beta$ exists, then it is unique. Indeed, any other such isomorphism would take the form $\gamma \beta$ for some global section $\gamma$ of $\Ow_X^{\times}$ such that $i^* \gamma = 1$. Since $f_* \Ow_X = \Ow_S$ (cf. \ref{chap2coh}), we have $\gamma = f^* \delta$ for some global section $\delta$ of $\Ow_S^{\times}$. Since the restriction of $\delta$ along the finite flat surjective morphism $Y \rightarrow S$ is trivial, one must have $\delta = 1$ as well, hence $\gamma = 1$.

\begin{prop}\label{chap2rigidrepr} Let $\Pic_{S}(X,Y)$ be the presheaf of abelian groups on $\Sch_{/S}^{\fp}$ which maps a finitely presented $S$-scheme $T$ to the set of isomorphism classes of $Y_T$-rigidified line bundles on $X_T$. Then, the presheaf $\Pic_{S}(X,Y)$ is representable by a smooth separated $S$-scheme of relative dimension $N +g -1$.
\end{prop}


We first consider the case where $N = 1$:

\begin{lem}\label{chap2rigidrepr2} The conclusion of Proposition \ref{chap2rigidrepr} holds if $N = 1$.
\end{lem}

Indeed, if $N=1$ then $Y$ is the image of a section $x : S \rightarrow X$ of $f$.  For any finitely presented $S$-scheme $T$, we have a morphism
\begin{align*}
\Pic(X \times_S T) &\rightarrow \Pic_{S}(X,x)(T) \\
\Lc &\rightarrow (\Lc \otimes (f^*x^* \Lc)^{-1}, \mathrm{id}).
\end{align*}
The kernel of this homomorphism consists of all invertible $\Ow_{X \times_S T}$-modules which are given by the pullback of an invertible $\Ow_T$-module. Moreover, any isomorphism class $(\Lc,\alpha)$ in $\Pic_{S}(X,x)(T)$ is the image of $\Lc$ by this morphism, hence its surjectivity. We conclude by \ref{chap2exactseq2} that the canonical projection morphism 
\begin{align*}
\Pic_{S}(X,x) &\rightarrow \Pic_{S}(X)\\
(\Lc,\alpha) &\rightarrow  \Lc,
\end{align*} 
is an isomorphism of presheaves of abelian groups on $\Sch_{/S}^{\fp}$. This yields Lemma \ref{chap2rigidrepr2} since $ \Pic_{S}(X)$ is a smooth separated $S$-scheme of relative dimension $g$ (cf. \ref{chap21.0.1}).

We now prove Proposition \ref{chap2rigidrepr}. Since $X \times_S Y \rightarrow Y$ has a section $x = (i \times \id_Y) \circ \Delta_Y$ where $\Delta_Y : Y \rightarrow Y \times_S Y$ is the diagonal morphism of $Y$, we deduce from Lemma \ref{chap2rigidrepr2} and its proof that the canonical projection morphism
$$
\Pic_{Y}(X \times_S Y,x) \rightarrow \Pic_{Y}(X \times_S Y) = \Pic_{S}(X) \times_S Y
$$
sending a pair $(\Lc, \alpha)$ to the class of $\Lc$ is an isomorphism. Let $Z$ be the $Y$-scheme $\Pic_{Y}(X \times_S Y,x)$, and let $(\Lc_u, \alpha_u)$ be the universal $x$-rigidified line bundle on $X \times_S Z$. The morphism $Y \times_S Z \rightarrow Z$ is finite locally free of rank $N$, so that the pushforward $\A$ (resp. $\M$) of $ \Ow_{Y \times_S Z}$ (resp. $i_Z^* \Lc_u$) is a locally free $\Ow_Z$-algebra of rank $N$ (resp. a locally free $\Ow_Z$-module of rank $N$). Let $\lambda : \M \rightarrow \Ow_Z$ be the surjective $\Ow_Z$-linear homomorphism corresponding to $\alpha_u^{-1} : x_Z^* \Lc_u \rightarrow \Ow_{Z}$.

Let $T$ be a $Y$-scheme, and let $(\Lc,\beta)$ be a $Y_T$-rigidified line bundle on $X_T$. The section $x_T : T \rightarrow X_T$ uniquely factors through $Y_T$ and we still denote by $x_T$ the corresponding section of $Y_T$. The pair $(\Lc,x_T^{*} \beta)$ is then an $x_T$-rigidified line bundle on $X_T$, so that there is a unique morphism $z : T \rightarrow Z$ such that $(\Lc,x_T^{*} \beta)$ is equivalent to the pullback by $z$ of $(\Lc_u, \alpha_u)$. Let us assume that $(\Lc,x_T^{*} \beta)$ is equal to this pullback. Then the global section $\beta$ of $i_T^* \Lc$ over $Y\times_S T$ provides a global section of $z^* \M$ over $T$, which we still denote by $\beta$, such that $(z^*\lambda)(\beta) = 1$ and $z^* \M = (z^* \A) \beta $. Conversely, any such section produces a $Y_T$-rigidification of $\Lc$ on $X_T$. The functor $\Pic_{S}(X,Y) \times_S Y = \Pic_{Y}(X \times_S Y, Y \times_S Y)$ is therefore isomorphic to the functor
\begin{align*}
\Sch_{/S}^{\fp} &\rightarrow \mathrm{Sets} \\
T &\mapsto \{ (z, \beta) \ | \ z \in Z(T), \beta \in \Gamma(T,z^* \M), \ \lambda(\beta)=1 \text{ and } \M_T = \A_T \beta \}.
\end{align*}
This implies that $\Pic_{S}(X,Y) \times_S Y$ is representable by a relatively affine $Z$-scheme, smooth of relative dimension $N-1$ over $Z$. By $\fppf$-descent of affine morphisms of schemes along the $\fppf$-cover $\Pic_{S}(X) \times_S Y \rightarrow \Pic_{S}(X)$, this implies the representability of $\Pic_{S}(X,Y)$ by an $S$-scheme, which is relatively affine and smooth of relative dimension $N-1$ over $\Pic_{S}(X)$. Since $\Pic_{S}(X)$ is separated and smooth of relative dimension $g$ over $S$ (cf. \ref{chap21.0}), the $S$-scheme $\Pic_{S}(X,Y)$ is separated and smooth of relative dimension $g+ N-1$.

\subsection{\label{chap21.2}} Let $f : X \rightarrow S$ be as in \ref{chap21.0}, and let $i : Y \hookrightarrow X$ be a closed subscheme of $X$, which is finite locally free over $S$ of degree $N \geq 1$, and let $U = X \setminus Y$ be its complement. A \textbf{$Y$-trivial effective Cartier divisor of degree $d$ on $X$} is a pair $(\Lc,\sigma)$ such that $\Lc$ is a locally free $\Ow_X$-module of rank $1$ and $\sigma : \Ow_X \hookrightarrow \Lc$ is an injective homomorphism such that $i^* \sigma$ is an isomorphism and such that the closed subscheme $V(\sigma)$ of $X$ defined by the vanishing of the ideal $\sigma \Lc^{-1}$ of $\Ow_X$ is finite locally free of rank $d$ over $S$. Two $Y$-trivial effective divisors $(\Lc,\sigma)$ and $(\Lc',\sigma')$ are \textbf{equivalent} if there is an isomorphism $\beta : \Lc \rightarrow \Lc'$ of $\Ow_X$-modules such that $\beta \sigma = \sigma'$. As in \ref{chap21.1}, if such an isomorphism exists then it is unique.

\begin{prop}\label{chap2symrep0}The map $(\Lc,\sigma) \mapsto (V(\sigma) \hookrightarrow X)$ is a bijection from the set of equivalence classes of $Y$-trivial effective Cartiers divisor of degree $d$ on $X$ onto the set of closed subschemes of $U$ which are finite locally free of degree $d$ over $S$.
\end{prop}

Let $(\Lc,\sigma)$ be a $Y$-trivial effective divisor of degree $d$ on $X$. The ideal $\mathcal{I} = \sigma \Lc^{-1}$ is an invertible ideal of $\Ow_{X}$ such that the vanishing locus $V(\mathcal{I})$ is finite locally free of rank $d$ over $S$ and is contained in $U$. The pair $(\Lc,\sigma)$ is equivalent to $(\mathcal{I}^{-1},1)$, and $\mathcal{I}$ is uniquely determined by $V(\mathcal{I})$. Conversely for any closed subscheme $Z$ of $U$ which is finite locally free of rank $d$ over $S$, the scheme $Z$ is proper over $S$ hence closed in $X$ as well, and its defining ideal $\mathcal{I}$ in $\Ow_{X_T}$ is invertible by (\cite{BLR}, 8.2.6(ii)). The pair $(\mathcal{I}^{-1},1)$ is then a $Y$-trivial effective Cartier divisor of degree $d$ on $X$.

\begin{prop}\label{chap2symrep} Let $d$ be an integer and let $\Div^{d,+}_{S}(X,Y)$ be the functor which to an $S$-scheme $T$ associates the set of equivalence classes of $Y_T$-trivial effective Cartier divisors of degree $d$ on $X_T$. Then $\Div^{d,+}_{S}(X,Y)$ is representable by the $S$-scheme $\Sy_S^d(U)$, the $d$-th symmetric power of $U = X \setminus Y$ over $S$ (cf. \ref{chap23.4}). In particular $\Div^{d,+}_{S}(X,Y)$ is smooth of relative dimension $d$ over $S$.
\end{prop}

By Proposition \ref{chap2symrep0}, the functor $\Div^{d,+}_{S}(X,Y)$ is isomorphic to the functor which sends an $S$-scheme $T$ to the set of closed subschemes of $U_T$ which are finite locally free of rank $d$ over $T$. In other words, $\Div^{d,+}_{S}(X,Y)$ is isomorphic to the Hilbert functor of $d$-points in the $S$-scheme $U$.

If $x$ is a $T$-point of $U$, we denote by $\Ow(-x)$ the kernel of the homomorphism $\Ow_{X \times_S T} \rightarrow x_* \Ow_{T}$, which is an invertible ideal sheaf, and by $\Ow(x)$ its dual, which is endowed with a section $1_x : \Ow_T \hookrightarrow \Ow(x)$. The morphism
\begin{align*}
\Sy_S^d(U) &\rightarrow \Div^{d,+}_{S}(X,Y) \\
(x_1,\dots,x_d) &\rightarrow \left( \bigotimes_{i=1}^d \Ow(x_i), \prod_{i=1}^d 1_{x_i} \right)
\end{align*}
is then an isomorphism of $\fppf$-sheaves by (\cite{SGA4}, XVII.6.3.9), hence Proposition \ref{chap2symrep}.

\begin{rema}\label{chap2compatib} Let $T$ be an $S$-scheme. Let $Z$ be a closed subscheme of $U_T$ which is finite locally free of rank $d$ over $T$, therefore defining a $T$-point of $\Div^{d,+}_{S}(X,Y) = \Sy_S^d(U)$ by Proposition \ref{chap2symrep0}. By (\cite{SGA4}, XVII.6.3.9), this $T$-point is given by the composition
$$
T \rightarrow \Sy_T^d(Z) \rightarrow \Sy_T^d(U_T) \rightarrow \Sy_S^d(U),
$$
where the first morphism is the canonical morphism from Proposition \ref{chap2morphalg}.
\end{rema}

\begin{prop}\label{chap2abeljacobi}
Let $d \geq N + 2g -1$ be an integer, and let $\Pic^d_{S}(X,Y)$ be the inverse image of $\Pic^d_{S}(X)$ by the natural morphism $\Pic_{S}(X,Y) \rightarrow \Pic_{S}(X)$. Then the Abel-Jacobi morphism
\begin{align*}
\Phi_d : \Div^{d,+}_{S}(X,Y) &\rightarrow \Pic^d_{S}(X,Y) \\
(\Lc,\sigma) &\mapsto (\Lc,i^*\sigma)
\end{align*}
is surjective smooth of relative dimension $d - N - g +1$ and it has geometrically connected fibers.
\end{prop}

Let $Z$ be the scheme $\Pic^d_{S}(X,Y)$, and let $(\Lc_u,\alpha_u)$ be the universal $Y$-rigidified line bundle of degree $d$ on $X_Z$. By (\cite{BLR}, 8.2.6(ii)), the closed subscheme $Y_Z$ of $X_Z$ is defined by an invertible ideal sheaf $\mathcal{I}$.

Let $\E$ be the pushforward of $\M = \Lc_u \otimes_{\Ow_{X_Z}} \mathcal{I}$ by the morphism $f_Z : X_Z \rightarrow Z$. By \ref{chap2coh2}, the $\Ow_Z$-module $\E$ is locally free of rank $d-N-g+1$, and for any morphism $T \rightarrow Z$ the canonical homomorphism
\begin{align*}
\E \otimes_{\Ow_Z} \Ow_T \rightarrow f_{T*} (\M \otimes_{\Ow_{X_Z}} \Ow_{X_T} )
\end{align*}
is an isomorphism, where $f_T : X_T \rightarrow T$ is the base change of $f$ by the morphism $T \rightarrow S$. We thus obtain an isomorphism
\begin{align}\label{chap2kl1}
 E  \rightarrow E',
\end{align}
of functors on the category of $Z$-schemes, where $E$ is the functor $T \mapsto \Gamma(T,\E \otimes_{\Ow_Z} \Ow_T)$ and $E'$ is the functor $T \mapsto \Gamma(X_T,\M \otimes_{\Ow_{X_Z}} \Ow_{X_T})$.
Let $\F$ be the pushfoward of $\Lc_u$ by the morphism $f_Z$. By the same argument, we obtain that the $\Ow_Z$-module $\F$ is locally free of rank $d-g+1$, and that we have an isomorphism
\begin{align}\label{chap2kl2}
 F  \rightarrow F',
\end{align}
of functors on the category of $Z$-schemes, where $F$ is the functor $T \mapsto \Gamma(T,\F \otimes_{\Ow_Z} \Ow_T)$ and $F'$ is the functor $T \mapsto \Gamma(X_T,\Lc_u \otimes_{\Ow_{X_Z}} \Ow_{X_T})$. Let us consider the exact sequence
$$
0 \rightarrow \M \rightarrow \Lc_u \rightarrow \Lc_u \otimes_{\Ow_{X_Z}} \Ow_{Y_Z} \rightarrow 0.
$$
Since $R^1f_{Z*} \M = 0 $ by \ref{chap2coh2}, we obtain an exact sequence
$$
0 \rightarrow \E \rightarrow \F \xrightarrow[]{} \mathcal{G} \rightarrow 0,
$$
where $\G$ is a locally free $\Ow_Z$-module of rank $N$. Together with \eqref{chap2kl1} and \eqref{chap2kl2}, this yields an exact sequence
$$
0 \rightarrow E' \rightarrow F' \xrightarrow[]{b} G \rightarrow 0,
$$
of $Z$-group schemes in $Z_{\fppf}$, where $G$ is the functor $T \mapsto \Gamma(T,\G_T \otimes_{\Ow_Z} \Ow_T)$. The section $\alpha_u$ of $\G$ over $Z$ corresponds to a morphism $\alpha_u : Z \rightarrow G$, and we have a morphism
\begin{align*}
\Div^{d,+}_{S}(X,Y) &\rightarrow F' \times_{b,G,\alpha_u} Z \\
(\Lc,\sigma) &\mapsto (\sigma, (\Lc,i^*\sigma)),
\end{align*}
which is an isomorphism: indeed, if $(\sigma, (\Lc, i^* \sigma))$ is a $T$-point of $F' \times_{b,G,\alpha_u} Z$, then for any point
$t$ of $T$ the restriction $\sigma_t$ of $\sigma$ to the fiber $X_t = X_T \times_T t$ is a global section of the line bundle
$\Lc_t = \Lc \otimes_{\Ow_{X_T}} \Ow_{X_t}$, which is non zero since non vanishing on $Y_t = Y_T \times_T t$, so that $\sigma_t : \Ow_{X_t} \rightarrow \Lc_t$ is an injective homomorphism and (\cite{EGA42}, 11.3.7) ensures that $\sigma : \Ow_X \rightarrow \Lc$ is an effective
Cartier divisor on the relative curve $X_T$, see also (\cite{BLR}, 8.2.6(iii)). Since $b$ is an $E'$-torsor over $G$ in $Z_{\fppf}$, we obtain that $\Div^{d,+}_{S}(X,Y)$ is an $E'$-torsor in $Z_{\fppf}$. Since $E'$ is isomorphic to $E$ by \eqref{chap2kl1}, it is smooth of relative dimension $d-N-g+1$ over $Z$ with geometrically connected fibers, hence the conclusion of Proposition \ref{chap2abeljacobi}.

\section{Geometric global class Field Theory}

\subsection{\label{chap22.2}}  Let $f : X \rightarrow S$ be a smooth morphism of schemes of relative dimension $1$, with connected geometric fibers of genus $g$, which is Zariski-locally projective over $S$, and let $i : Y \hookrightarrow X$ be a closed subscheme of $X$ which is finite locally free over $S$ of degree $N \geq 1$. Let $j : U \rightarrow X$ be the open complement of $Y$. Let $\Lambda$ be a finite ring whose cardinality is invertible on $S$.

\begin{defi}\label{chap2ramif6} A locally free $\Lambda$-module $\F$ of rank $1$ in $U_{\et}$ has \textbf{ramification bounded by $Y$ over $S$} if for any geometric point $\bar{x}$ of $Y$ with image $\bar{s}$ in $S$, the restriction of $\F$ to $\Spec(\widehat{\Ow_{X_{\bar{s}},\bar{x}}}) \times_{X_{\bar{s}}} U_{\bar{s}}$ has ramification bounded by the multiplicity of $Y_{\bar{s}}$ at $\bar{x}$ (cf. \ref{chap2ramif}).
\end{defi}

\begin{teo}\label{chap2GCFT} Let $\F$ be a locally free $\Lambda$-module of rank $1$ in $U_{\et}$ with ramification bounded by $Y$ over $S$ (cf. \ref{chap2ramif6}). Then, there is a unique (up to isomorphism) multiplicative locally free $\Lambda$-module $\G$ of rank $1$ on the $S$-group scheme $\Pic_{S}(X,Y)$ (cf. \ref{chap2charsheaf}) such that the pullback of $\G$ by the Abel-Jacobi morphism
$$
U \rightarrow \Pic_{S}(X,Y),
$$
which sends $x$ to $(\Ow(x), 1)$, is isomorphic to $\F$.
\end{teo}

In Section \ref{chap22.0}, we study the restriction of the locally free $\Lambda$-module $\F^{[d]}$ of rank $1$ on $\Div^{d,+}_{S}(X,Y)$ (cf. \ref{chap2descent} and \ref{chap2symrep}) to a geometric fiber of the Abel-Jacobi morphism (cf. \ref{chap2abeljacobi})
\begin{align*}
\Phi_d : \Div^{d,+}_{ S}(X,Y) &\rightarrow \Pic_{S}^d(X,Y) \\
(\Lc,\sigma) &\mapsto (\Lc,i^*\sigma).
\end{align*}
This study will enable us to prove Theorem \ref{chap2GCFT} in Section \ref{chap22.3}.

\subsection{\label{chap22.0}} Let $k$ be an algebraically closed field, let $X$ be a smooth connected projective curve of genus $g$ over $k$ and let $i : Y \rightarrow X$ be an effective Cartier divisor of degree $N$ with complement$U$ in $X$. Let $\Lc$ be a line bundle of degree $d \geq N + 2g-1$ on $X$, and let $V$ be the $(d-N-g+1)$-dimensional affine space over $k$ associated to the $k$-vector space $\Vc = H^0(X,\Lc(-Y))$, i.e. $V$ is the spectrum of the symmetric algebra of the $k$-module $\mathrm{Hom}_k(\Vc,k)$. Let $\tau$ be a global section of $\Lc$ on $X$ such that $i^* \tau : \Ow_Y \rightarrow i^* \Lc$ is an isomorphism.
%

\begin{prop}\label{chap2ramif4} Let $\Lambda$ be a finite ring of cardinality invertible in $k$, and let $\F$ be a locally free $\Lambda$-module of rank $1$ in $U_{\et}$, with ramification bounded by $Y$ (cf. \ref{chap2ramif6}). Then the pullback of $\F^{[d]}$ (cf. \ref{chap2descent}) by the morphism
$$
V \rightarrow \Div^{d,+}_{k}(X,Y),
$$
which sends a section $s$ of $V$ to $(\Lc, \tau - s)$, is a constant \'etale sheaf.
\end{prop}

The morphism
$$
V \rightarrow \Div^{d,+}_{k}(X,Y),
$$
which sends a point $\sigma$ of $V$ to $(\Lc, \tau - \sigma)$, is an isomorphism from $V$ to the fiber of $\Phi_d$ over the $k$-point $(\Lc,i^* \tau)$, cf. \ref{chap2abeljacobi}. Proposition \ref{chap2ramif4} thus implies:

\begin{cor}\label{chap2ramif15} Let $\F$ be as in Proposition \ref{chap2ramif4}. Then the locally free $\Lambda$-module $\F^{[d]}$ on $\Div^{d,+}_{k}(X,Y)_{\et}$ is constant on the fiber at $(\Lc,i^* \tau)$ of the morphism
$$
\Phi_d : \Div^{d,+}_{k}(X,Y) \rightarrow \Pic^{d}_{k}(X,Y)
$$
from \ref{chap2abeljacobi}.
\end{cor}

We now prove Proposition \ref{chap2ramif4}. To this end, we consider the morphism
$$
\psi : \mathbb{A}^1_V \rightarrow \Div^{d,+}_{k}(X,Y),
$$
which sends a pair $(t,\sigma)$, where $t$ and $\sigma$ are points of $\mathbb{A}^1_k$ and $V$ respectively, to the point $(\Lc, \tau - t \sigma)$ of $\Div^{d,+}_{k}(X,Y)$. Let $\F$ be as in Propoosition \ref{chap2ramif4}, and let $\G$ be the pullback by $\psi$ of $\F^{[d]}$ (cf. \ref{chap2descent}). Denoting by $\iota_t : V \rightarrow \mathbb{A}^1_V$ the section corresponding to an element $t$ of $k = \mathbb{A}_k^1(k)$, we must prove that the sheaf $\iota_1^{-1} \G$ is constant. The sheaf $\iota_0^{-1} \G$ is constant, since $\psi \iota_0$ is a constant morphism, hence it is sufficent to prove that $\iota_1^{-1} \G$ and $\iota_0^{-1} \G$ are isomorphic. The latter fact follows from the following lemma:

\begin{lem}\label{chap2tame} The locally free $\Lambda$-module $\G$ is the pullback of an etale sheaf on $V$ by the projection $\pi : \mathbb{A}^1_V \rightarrow V$.
\end{lem}

We now prove Lemma \ref{chap2tame}. We start by proving that $\G$ is constant on each geometric fiber of the projection $\pi$. Since the formation of $\psi$ and $\G$ is compatible with the base change along any field extension of $k$, it is sufficient to show that $\G$ is constant on each fiber of the projection $\mathbb{A}^1_V \rightarrow V$ at a $k$-point $\sigma$ of $V$. If $\sigma = 0$, then the restriction of $\psi$ to the fiber of $\pi$ above $\sigma$ is constant, hence $\G$ is constant on this fiber.

We now assume that $\sigma$ is non zero. Since $\sigma$ vanishes on the non empty divisor $Y$ and $\tau$ does not, the sections $\sigma$ and $\tau$ are $k$-linearly independent in $H^0(X,\Lc)$. Let $D$ be the greatest divisor on $X$ such that $D \leq \mathrm{div}(\sigma)$ and $D \leq \mathrm{div}(\tau)$. Since the divisor of $\tau$ is contained in $U$, so is $D$. We can then write $\sigma = \widetilde{\sigma} 1_D$ and $\tau = \widetilde{\tau} 1_D$, where $1_D$ is the canonical section of $\Ow(D)$ and $\widetilde{\sigma} , \widetilde{\tau}$ are global sections of $\Lc(-D)$ on $X$ without common zeroes. Thus $f = [\widetilde{\tau} : \widetilde{\sigma} ]$ is a well defined non constant morphism from $X$ to $\mathbb{P}_{k}^1$. Thus, if $W$ is the closed subscheme of $X \times_k \mathbb{A}^1_k$ defined by the vanishing of $\tau - t \sigma$, where $t$ is the coordinate on $\mathbb{A}^1_k$, then we have
$$
W = D \times_{k} \mathbb{A}^1_k \cup (\mathrm{Graph}(f) \cap X \times_k \mathbb{A}^1_k)  \hookrightarrow U \times_{k} \mathbb{A}^1_k.
$$
Moreover, the projection $W \rightarrow \mathbb{A}^1_k$ is finite flat of degree $d$, and the restriction of $\psi$ to the fiber at $\sigma$ factors as
$$
 \mathbb{A}^1_k \xrightarrow[]{\varphi} \Sym_{\mathbb{A}^1_k}^d( W) \rightarrow \Sym_{\mathbb{A}^1_k}^d(U \times_{k} \mathbb{A}^1_k) \rightarrow  \Sym_{k}^d(U) \rightarrow \Div^{d,+}_{k}(X,Y),
$$
where the first morphism $\varphi$ is obtained from Proposition \ref{chap2morphalg}, and the last morphism is the isomorphism from Proposition \ref{chap2symrep}. Moreover, the pullback of $\F^{[d]}$ to $ \Sym_{\mathbb{A}^1_k}^d( W)$ coincides with $(p_1^{-1} \F)^{[d]}$, where $p_1 : W \rightarrow U$ is the first projection. In particular, the sheaf $\G$ is isomorphic to $\varphi^{-1} (p_1^{-1} \F)^{[d]}$.

Let $K=  k((t^{-1}))$ and let $\eta = \Spec(K) \rightarrow \mathbb{A}^1_{k}$ be the corresponding punctured formal neighbourhood of $\infty$. Let us form the following commutative diagram:

\begin{center}
 \begin{tikzpicture}[scale=1.5]

\node (A) at (0,0) {$\eta$};
\node (B) at (0,1) {$\mathbb{A}^1_{k}$};
\node (C) at (2,0) {$\Sym^d_{\eta}(W \times_{\mathbb{A}^1_{k}} \eta)$.};
\node (D) at (2,1) {$\Sym^d_{\mathbb{A}^1_{k}}(W )$};
\path[->,font=\scriptsize]
(A) edge (B)
(B) edge node[above]{$\varphi$} (D)
(C) edge (D) 
(A) edge (C) ;
\end{tikzpicture} 
\end{center}

We can then write
$$
W \times_{\mathbb{A}^1_{k}} \eta = D \times_{k} \eta \cup  \mathrm{Graph}(f) \times_{\mathbb{P}_k^1} \eta = D \times_{k} \eta \cup X \times_{f,\mathbb{P}_{k}^1} \eta.
$$
The divisors $ D \times_{k} \eta$ and $X \times_{f,\mathbb{P}_{k}^1} \eta$ of $X \times_{k} \eta$ are disjoint, since the former lies over closed points of $X$, while the latter lies over the generic point of $X$. We thus have a decomposition
$$
W \times_{\mathbb{A}^1_{k}} \eta = D \times_{k} \eta \amalg  X \times_{f,\mathbb{P}_{k}^1}\eta = \coprod_{i} \Spec(L_i)
$$
where $L_i$ is either of the form $K[T]/(T^{d_i})$ if $\Spec(L_i)$ is a connected component of $D \times_{k} \eta$, or a field extension of degree $d_i$ of $K$ if $\Spec(L_i)$ is a connected component of $X \times_{f,\mathbb{P}_{k}^1}\eta$. In the former case, the restriction of $p_1^{-1} \F$ to $\Spec(L_i)$ is constant, while in the latter case, we have the further information that the restriction of $p_1^{-1} \F$ to $\Spec(L_i)$ has ramification bounded by $d_i$ (cf. \ref{chap2ramif}), since the ramification index of $f$ at a point $x$ above $\infty$ is greater than or equal to the multiplicity of $Y$ at $x$, and $\F$ has ramification bounded by $Y$ by assumption. Moreover, we have $\sum_{i} d_i = d$, and the morphism $\eta \rightarrow\Sym^d_{\eta}(W \times_{\mathbb{A}^1_{k}} \eta)$ factors through the canonical morphism
$$
\prod_i \Sym^{d_i}_{\eta}(\Spec(L_i)) \rightarrow \Sym^d_{\eta}(W \times_{\mathbb{A}^1_{k}} \eta).
$$ 	
By \ref{chap2ramif2}, we obtain that the restriction of $\G$ to $\eta$ is tamely ramified. Since the tame fundamental group of $\mathbb{A}^1_{k}$ is trivial, we conclude that $\G$ is a constant \'etale $\Lambda$-module on the fiber of $\pi$ at $\sigma$. The conclusion of Lemma \ref{chap2tame} then follows from a descent result, namely Lemma \ref{chap2pull} below.

\begin{rema} While the proof of Proposition \ref{chap2ramif2}, which constitutes the core of the proof of Lemma \ref{chap2tame} above, uses geometric local class field theory, it should be noticed that its statement does not refer to it. This explains why no form of local-global compatibility is required in the proof of Lemma \ref{chap2tame}.
\end{rema}

%
%
%
%

\begin{lem}\label{chap2pull} Let $g : T' \rightarrow T$ be a quasi-compact smooth compactifiable morphism of schemes of relative dimension $\delta$ with geometrically connected fibers, and let $\G$ be an \'etale sheaf of $\Lambda$-modules on $T'_{\eet}$ which is constant on each geometric fiber of $g$. Then $\G$ is isomorphic to the pullback by $g$ of an \'etale sheaf of $\Lambda$-modules on $T_{\eet}$. 
\end{lem}

By (\cite{SGA4}, XVIII 3.2.5) the functor $Rg_{!}$ on the derived category of $\Lambda$-modules on $T$ admits the functor $g^{!} : K \mapsto g^* K(\delta)[2 \delta]$ as a right adjoint. Let us apply the functor $\Hm^0$ to the adjunction morphism $\G \rightarrow g^! Rg_! \G$. The morphism
$$
\G \rightarrow \Hm^0(g^! Rg_! \G)  = g^* R^{2 \delta} g_! \G (\delta)
$$
is an isomorphism, as can be seen by checking the stalks at geometric points with the proper base change theorem.


\subsection{\label{chap22.3}} We now prove Theorem \ref{chap2GCFT}. Let $\F$ be a locally free $\Lambda$-module of rank $1$ over $U_{\et}$. The family $(F^{[d]})_{d \geq 0}$ of locally free $\Lambda$-modules of rank $1$ yields a multiplicative \'etale $\Lambda$-module of rank $1$ over the $S$-semigroup scheme
$$
\Div^{+}_{S}(X,Y) = \coprod_{d \geq 0} \Div^{d,+}_{ S}(X,Y).
$$
For each integer $d \geq N +2g-1$, Corollary \ref{chap2ramif15} implies that the locally free $\Lambda$-module $\F^{[d]}$ of rank $1$ on $\Div^{d,+}_{S}(X,Y)$ (cf. \ref{chap2descent} and \ref{chap2symrep}) is constant on the geometric fibers of the smooth surjective morphism (cf. \ref{chap2abeljacobi})
\begin{align*}
\Phi_d : \Div^{d,+}_{ S}(X,Y) &\rightarrow \Pic_{S}^d(X,Y) \\
(\Lc,\sigma) &\mapsto (\Lc,i^*\sigma).
\end{align*}
This morphism satisfies the conditions of Lemma \ref{chap2pull} by Proposition \ref{chap2abeljacobi}. We can therefore apply Lemma \ref{chap2pull}, and we obtain a locally free $\Lambda$-module $\G_d$ of rank $1$ over $\Pic_{S}^d(X,Y)$ such that $\Phi_d^{-1} \G_d$ is isomorphic to $\F^{[d]}$. By Proposition \ref{chap2mondesc}, the family $(\G_d)_{d \geq N+2g-1}$ yields a multiplicative locally free $\Lambda$-module of rank $1$ on the $S$-semigroup scheme
$$
M = \coprod_{d \geq N+2g-1} \Pic_{S}^d(X,Y).
$$
Since the morphism 
\begin{align*}
\rho : M \times_S M &\rightarrow \Pic_{S}(X,Y) \\
(x,y) &\mapsto x y^{-1}
\end{align*}
is faithfully flat and quasi-compact, we can apply Proposition \ref{chap2extension}, which yields a multiplicative locally free $\Lambda$-module $\G$ of rank $1$ over $\Pic_{S}(X,Y)$ whose restriction to $\Pic_{S}^d(X,Y)$ coincides with $\G_d$ for $d \geq N + 2g-1$. The families $(\F^{[d]})_{d \geq 0}$ and $(\Phi_d^{-1} \G_d)_{d \geq 0}$ yield multiplicative locally free $\Lambda$-modules of rank $1$ on the $S$-semigroup scheme $\Div^{+}_{S}(X,Y) = \coprod_{d \geq 0} \Div^{d,+}_{ S}(X,Y)$, whose restrictions to the ideal 
$$
I = \coprod_{d \geq N+2g-1} \Div^{d,+}_{ S}(X,Y)
$$
of $\Div^{+}_{S}(X,Y)$ are isomorphic. We obtain by Proposition \ref{chap2mondesc2} an isomorphism from $\F^{[d]}$ to $\Phi_d^{-1} \G_d$ for each $d \geq 0$. In particular, the locally free $\Lambda$-module $\Phi_1^{-1} \G_1$ of rank $1$ is isomorphic to $\F$.

\bibliographystyle{amsalpha}

\begin{thebibliography}{2}
%

\bibitem[BE01]{bloch}
S. Bloch, H. Esnault, ``Gauss-Manin determinants for rank 1 irregular connections on curves'',
Mathematische Annalen, Volume 321, 15-87, 2001. With an addendum: the letter of P. Deligne to J.-P. Serre (Feb. 74) on $\varepsilon$-factors, 65-87. 
%


\bibitem[BLR90]{BLR} S. Bosch, W. Lutkebohmert and M. Raynaud, ``Néron Models'', Ergebnisse der Mathematik und ihrer
Grenzgebiete, Springer-Verlag, 1990.

\bibitem[CC13]{CC} C. Contou-Carr\`ere, Jacobienne locale d'une courbe formelle relative, Rendiconti del Seminario Matematico della Universit\`a di Padova 130 (2013), pp.1-106.

\bibitem[EGA3]{EGA3} A. Grothendieck, J. Dieudonné.  \emph{\'{E}l\'ements de g\'eom\'etrie alg\'ebrique. {III}. \'Etude cohomologique des faisceaux coh\'erents}, Publications Mathématiques de l'I.H.\'E.S 11-17, 1961-1963.

\bibitem[EGA4]{EGA42}
A. Grothendieck, \emph{\'{E}l\'ements de g\'eom\'etrie alg\'ebrique. {IV}.
  \'{E}tude locale des sch\'emas et des morphismes de sch\'emas.}, Publications Mathématiques de l'I.H.\'E.S 20-24-28-32,
  1964-1967.

\bibitem[SGA1]{SGA1}
A. Grothendieck, \emph{S\'eminaire de G\'eom\'etrie Alg\'ebrique du Bois Marie - 1960-61 - Rev\^etements \'etales et groupe fondamental - (SGA 1)}, Springer-Verlag, LNM 224, 1971.

\bibitem[SGA4]{SGA4}
A. Grothendieck, \emph{S\'eminaire de G\'eom\'etrie Alg\'ebrique du Bois Marie - 1963-64 - Th\'eorie des topos et cohomologie \'etale des sch\'emas - (SGA 4)}, Springer-Verlag, LNM 269/270/305, 1972/3. 

%



%
%
%
%
%
%
%



%

%
%
%

\bibitem[La56]{Lang} S. Lang, ``Sur les s\'eries $L$ d'une vari\'et\'e alg\'ebrique'', Bulletin de la S. M. F. 84 (1956), p. 385-407.

%

\bibitem[La90]{Laumon} G. Laumon, ``Faisceaux automorphes li\'es aux s\'eries d'Eisenstein'', in Automorphic forms, Shimura varieties, and L-functions, Vol. I, 227-281, Perspect. Math., 10, Academic Press (1990).

\bibitem[MB85]{MB} L. Moret-Bailly, ``Pinceaux de vari\'et\'es ab\'eliennes'', Ast\'erisque 129 (1985).


\bibitem[Ro54]{Ros} M. Rosenlicht, ``Generalized Jacobian Varieties'', Annals of Mathematics 59 (1954), p. 505-530.

\bibitem[Se59]{JPS3} J.-P. Serre, ``Groupes alg\'ebriques et corps de classes'', Hermann (Paris), 1959.

\bibitem[Se61]{JPS} J.-P. Serre, ``Sur les corps locaux \`a corps r\'esiduel alg\'ebriquement clos'', Bulletin de la S. M. F. 89 (1961), p. 105-154.

\bibitem[Se68]{JPS2} J.-P. Serre, ``Corps locaux'', Hermann (Paris), 1968.
%

%


\bibitem[Su13]{TS} T. Suzuki, ``Some remarks on the local class field theory of Serre and Hazewinkel'', Bulletin de la S.M.F. 141 (2013), p. 1-24.


\bibitem[SP]{Stacks} The Stacks project, https://stacks.math.columbia.edu, 2019.

\bibitem[Ta18]{T18} D. Takeuchi, ``Blow-ups and the class field theory for curves'', arxiv:1804.02136.

%
%
%
%
%
%
%
%
%
%
%
%
%
%
%
%
%
%

\end{thebibliography}

\end{document}